\newcount\ite\ite=1\def\0{\global\ite=1\1}
\def\1{\item{\rm(\romannumeral\the\ite)}\advance\ite1\quad}

\documentclass[12pt,oneside]{article}
\usepackage{amssymb}

\font\teneufm=eufm10 scaled \magstep1
\font\seveneufm=eufm7 scaled \magstep1
\font\fiveeufm=eufm5  scaled \magstep1
\newfam\eufmfam

\textfont\eufmfam=\teneufm
\scriptfont\eufmfam=\seveneufm
\scriptscriptfont\eufmfam=\fiveeufm

\newfam\msbfam
\font\tenmsb=msbm10 scaled \magstep1  \textfont\msbfam=\tenmsb
\font\sevenmsb=msbm7 scaled \magstep1 \scriptfont\msbfam=\sevenmsb
\font\fivemsb=msbm5 scaled \magstep1  \scriptscriptfont\msbfam=\fivemsb

\def\dd#1{\raise1.5pt\hbox{$\,\partial\!$}/\raise-2.5pt\hbox{$\!\partial#1\,$}}

\def\5#1{{\mathcal #1}}

\def\RR{{\mathbb R}}
\def\CC{{\mathbb C}}

\def\ZZ{{\mathbb Z}}

\def\PP{{\mathbb P}}

\def\ra{\rightarrow}

\def\GL{\mathop{\rm GL}\nolimits}
\def\Sp{\mathop{\rm Sp}\nolimits}
\def\SO{\mathop{\rm SO}\nolimits}

\def\mod{\mathop{\rm mod}\nolimits}

\def\aut{\mathop{\rm aut}\nolimits}
\def\hol{\mathop{\rm hol}\nolimits}
\def\Aut{\mathop{\rm Aut}\nolimits}
\def\aut{\mathop{\rm aut}\nolimits}
\def\Stab{\mathop{\rm Stab}\nolimits}
\def\Diff{\mathop{\rm Diff}\nolimits}
\def\Im{\mathop{\rm Im}\nolimits}
\def\Re{\mathop{\rm Re}\nolimits}

 \def\HollowBoxx #1#2#3{{\dimen0=#1 \advance\dimen0 by -#2
       \dimen1=#1 \advance\dimen1 by #3
        \vrule height 0pt depth #3 width #2
       \hskip -#3
       \vrule height #1 depth #3 width #3}}
 \def\LeftContraction{\mathord{\kern1.45pt \HollowBoxx{6pt}{3.5pt}{.4pt}}\,}

 \def\HollowBox #1#2#3{{\dimen0=#1 \advance\dimen0 by -#3
       \dimen1=#1 \advance\dimen1 by #3
        \vrule height #1 depth #3 width #3
        \vrule height 0pt depth #3 width #2
        \hskip -#3}}
 \def\RightContraction{\mathord{\, \HollowBox{6pt}{3.1pt}{.4pt}} \kern1.6pt}

\def\qed{{\hfill $\Box$}}
\newtheorem{theorem}{THEOREM}[section]
\newtheorem{corollary}[theorem]{Corollary}
\newtheorem{lemma}[theorem]{Lemma}

\newtheorem{remark}[theorem]{Remark}

\newtheorem{conjecture}[theorem]{Conjecture}

  \renewenvironment{thebibliography}[1]{%
    \begin{oldthebibliography}{#1}%
      \setlength{\parskip}{0.1ex}%
      \setlength{\itemsep}{0.5ex}%
  }%
  {%
    \end{oldthebibliography}%
}

\makeatletter
\def\blfootnote{\xdef\@thefnmark{}\@footnotetext}
\makeatother

\begin{document}

\begin{center}
{\Large \bf Reduction of Five-Dimensional\\
\vspace{0.1cm}
Uniformly Levi Degenerate CR Structures\\
\vspace{0.3cm}
to Absolute Parallelisms}
\medskip
\medskip\\
\normalsize Alexander Isaev\footnote{The research is supported by the Australian Research Council.} and Dmitri Zaitsev\footnote{Supported in part by Science Foundation Ireland grant 10/RFP/MTH2878.}
\end{center}

\begin{quotation} 
{\small \sl \noindent Let ${\mathfrak C}_{2,1}$ be the class of connected 5-dimensional CR-hypersurfaces that are 2-nondegenerate and uniformly Levi degenerate of rank 1. We show that the CR-structures in ${\mathfrak C}_{2,1}$ are reducible to ${\mathfrak{so}}(3,2)$-valued absolute parallelisms and give applications of this result.}
\end{quotation}

\thispagestyle{empty}

\pagestyle{myheadings}
\markboth{A. Isaev and D. Zaitsev}{Uniformly Levi Degenerate CR Structures}

\setcounter{section}{0}

\section{Introduction}\label{intro}
\setcounter{equation}{0}

An almost CR-structure on a smooth manifold $M$ is a subbundle $H(M)\subset T(M)$ of the tangent bundle of even rank endowed with operators of complex structure $J_p:H_p(M)\ra H_p(M)$, $J_p^2= -\hbox{id}$, that smoothly depend on $p$. A manifold equipped with an almost CR-structure is called an almost CR-manifold. The subspaces $H_p(M)$ are called the complex tangent spaces to $M$, and their complex dimension, denoted by $\hbox{CRdim}\, M$, is the CR-dimension of  $M$. The complementary dimension  $\hbox{CRcodim}\, M:=\dim M-2\hbox{CRdim}\, M$ is called the CR-codimension of $M$. Further, a smooth map $f:M\ra \tilde M$ between two almost CR-manifolds is a CR-map if for every $p\in M$ the differential $df(p)$ of $f$ at $p$ maps $H_p(M)$ into $H_{f(p)}(\tilde M)$ and is complex-linear on $H_p(M)$. If for two almost CR-manifolds $M$, $\tilde M$ of equal CR-dimensions there exists a diffeomorphism $f$ from $M$ onto $\tilde M$ that is also a CR map, then the manifolds are said to be CR-equivalent and $f$ is called a CR-isomorphism.

We are interested, in particular, in the equivalence problem for almost CR-manifolds. This problem can be viewed as a special case of the equivalence problem for ${\bf G}$-structures. Let ${\bf G}\subset \GL(d,\RR)$ be a Lie subgroup. A ${\bf G}$-structure on a $d$-dimensional manifold $M$ is a subbundle ${\mathcal S}$ of the frame bundle $F(M)$ over $M$ that is a principal ${\bf G}$-bundle. Two ${\bf G}$-structures ${\mathcal S}$, $\tilde {\mathcal S}$ on manifolds $M$, $\tilde M$, respectively, are called equivalent if there is a ${\bf G}$-structure isomorphism between ${\mathcal S}$, $\tilde{\mathcal S}$, i.e.~a diffeomorphism $f$ from $M$ onto $\tilde M$ such that the induced mapping\linebreak $f_{*}:F(M)\ra F(\tilde M)$ maps ${\mathcal S}$ onto $\tilde{\mathcal S}$. The almost CR-structure of a manifold $M$ of CR-dimension $n$ and CR-codimension $k$ is a ${\bf G}$-structure with ${\bf G}$ being the group of all nondegenerate linear transformations of $\CC^n\oplus\RR^k$ that preserve the first component and are complex-linear on it. The notion of equivalence of such ${\bf G}$-structures is then exactly that of almost CR-structures. For convenience, when speaking about ${\bf G}$-structures below, we replace the frame bundle $F(M)$ by the coframe bundle.

Important examples of ${\bf G}$-structures are $\{e\}$-structures, where $\{e\}$ is the one-element group. They are called absolute parallelisms, and on a  $d$-dimensional manifold $M$ any such structure is given by an $\RR^d$-valued 1-form that for every $p\in M$ defines an isomorphism between $T_p(M)$ and $\RR^d$. The (local) equivalence problem for absolute parallelisms is well-understood (see e.g.\linebreak p. 344 in \cite{Ste}), and therefore one may approach the equivalence problem for general\linebreak ${\bf G}$-structures by attempting to reduce them to absolute parallelisms. Let ${\mathfrak C}$ be a class of manifolds equipped with ${\bf G}$-structures. The ${\bf G}$-structures in ${\mathfrak C}$ are said to reduce to absolute parallelisms if to every $M\in{\mathfrak C}$ one can assign a fiber bundle ${\mathcal P}$ and an absolute parallelism $\pi$ on ${\mathcal P}$ in such a way that for any $M,\tilde M\in{\mathfrak C}$ the following holds: (i) any ${\bf G}$-structure isomorphism $f:M\ra \tilde M$ can be lifted to a diffeomorphism $F: {\mathcal P}\ra\tilde {\mathcal P}$ satisfying
\begin{equation}
F^{*}\tilde\pi=\pi,\label{eq8}
\end{equation}
and (ii) any diffeomorphism $F: {\mathcal P}\ra\tilde {\mathcal P}$ satisfying (\ref{eq8}) 
is a bundle isomorphism that is a lift of a ${\bf G}$-structure isomorphism $f:M\ra \tilde M$.
 
\'E. Cartan developed a general method for reducing ${\bf G}$-structures to absolute parallelisms, which works, in particular, for Riemannian and conformal structures (see \cite{C2}, \cite{IL}). However, in contrast with the conformal and Riemannian cases, Cartan's reduction procedure is not directly applicable to almost CR-structures. Nevertheless, as we will see below, in a number of situations almost CR-structures are known to reduce to parallelisms. In these cases the CR-structures are required to satisfy additional conditions, which we will now introduce. Broadly speaking, one needs assumptions of two kinds: integrability (or at least partial integrability) and some nondegeneracy.

We start with the integrability condition. Let $M$ be an almost CR-manifold. For every $p\in M$ consider the complexification
$H_p(M)\otimes_{\RR}\CC$ of the complex tangent space at $p$. It can be
represented as the direct sum
$$
H_p(M)\otimes_{\RR}\CC=H_p^{(1,0)}(M)\oplus H_p^{(0,1)}(M),
$$
where
$$
\begin{array}{l}
H_p^{(1,0)}(M):=\{X-iJ_pX:X\in H_p(M)\},\\
\vspace{-0.3cm}\\
H_p^{(0,1)}(M):=\{X+iJ_pX:X\in H_p(M)\}.
\end{array}
$$
Then the almost CR-structure on $M$ is said to be integrable if the bundle $H^{(1,0)}$ is involutive, i.e.~for any pair of local
sections ${\mathfrak z},{\mathfrak z}'$ of $H^{(1,0)}(M)$ the commutator $[{\mathfrak z},{\mathfrak z}']$
is also a local section of $H^{(1,0)}(M)$. An integrable almost CR-structure is called a CR-structure and a manifold equipped with a CR-structure a CR-manifold. 

Next, one requires some nondegeneracy condition. The most common condition of this kind arises from the Levi form, which comes from taking commutators of local sections of
$H^{(1,0)}(M)$ and $H^{(0,1)}(M)$. Let $p\in M$, $Z,Z'\in
H_p^{(1,0)}(M)$. Choose  local sections ${\mathfrak z}$, ${\mathfrak z}'$ of $H^{(1,0)}(M)$ near
$p$ such that ${\mathfrak z}(p)=Z$, ${\mathfrak z}'(p)=Z'$. The Levi form of $M$ at
$p$ is then the Hermitian form on $H_p^{(1,0)}(M)$ with values in $(T_p(M)/H_p(M))\otimes_{\RR}\CC$ given by
$$
{\mathcal L}_M(p)(Z,Z'):=i[{\mathfrak z},\overline{{\mathfrak z}'}](p)(\mod H_p(M)\otimes_{\RR}\CC).
$$
For fixed $Z$ and $Z'$ the right-hand side of the above formula is independent of the choice of ${\mathfrak z}$ and ${\mathfrak z}'$. The Levi form is often treated as a $\CC^k$-valued Hermitian form (i.e.~a vector of $k$ Hermitian forms), where $k:=\hbox{CRcodim}\,M$. As a $\CC^k$-valued Hermitian form, the Levi form is defined uniquely up to the choice of coordinates in $T_p(M)/H_p(M)$. Now, a CR-manifold $M$ is called Levi nondegenerate if its Levi form at any $p\in M$ is nondegenerate, where a $\CC^k$-valued Hermitian form $h=(h_1,\dots,h_k)$ on a vector space $V$ is said to be nondegenerate provided the following holds: (i) the scalar-valued forms $h_1,\dots,h_k$ are linearly independent over $\RR$, and\linebreak (ii) $h(v,v')=0$ for all $v'\in V$ implies $v=0$. 

The problem of reducing CR-structures to absolute parallelisms is best-studied for CR-hypersurfaces, i.e.~CR-manifolds of CR-codimension 1. The first result of this type goes back to \'E. Cartan (see \cite{C1}) who showed that reduction takes place for all 3-dimensional Levi nondegenerate CR-hyper\-surfaces (note that the method of \cite{C1} differs from Cartan's approach to general ${\bf G}$-structures mentioned earlier). Further, in \cite{T2} Tanaka proposed a reduction procedure that applies to many geometric structures, in particular to all Levi nondegenerate CR-hypersurfaces. However, Tanaka's work became widely known only after the publication of Chern-Moser's article \cite{CM} where the problem was solved independently in the Levi nondegenerate CR-hypersurface case (cf. \cite{T3}). Next, in a certain more general situation (namely for Levi nondegenerate partially integrable almost CR-structures of CR-codimension 1) reduction to parallelisms was obtained in \cite{CSc} as part of a general parabolic geometry approach (for more details on this approach see \cite{CSl}). Finally, there is a number of classes of Levi nondegenerate CR-manifolds of CR-codimension greater than 1 for which reduction to parallelisms has also been established. Indeed, the method of \cite{T2} yields reduction for strongly uniform CR-structures, where strong uniformity means that the Levi forms at all points are pairwise equivalent. Explicit constructions for special strongly uniform CR-structures were also given in \cite{EIS}, \cite{ScSl}, \cite{ScSp}. In the absence of strong uniformity reduction was obtained for some classes of CR-manifolds in \cite{GM}, \cite{L}, \cite{Mi}. We stress that all of the above results assume Levi nondegeneracy.

In this paper we set out to obtain reduction to absolute parallelisms for CR-hypersurfaces satisfying a different nondegeneracy condition, the so-called 2-nondegeneracy. For the general notion of $k$-nondegeneracy (as well as other nondegeneracy conditions) we refer the reader to Chapter XI in \cite{BER} and note that for CR-hypersurfaces 1-nondegeneracy is equivalent to Levi nondegeneracy. We consider 5-dimen\-sional CR-hypersurfaces that are uniformly Levi degenerate of rank 1. In the terminology of \cite{E} this means that the kernel $\ker{\mathcal L}_M(p)$ of the Levi form has dimension 1 at every $p\in M$, where 
$$
\ker{\mathcal L}_M(p):=\left\{Z\in H_p^{(1,0)}(M): {\mathcal L}_M(p)(Z,Z')=0\,\,\,\hbox{for all $Z'\in H_p^{(1,0)}(M)$}\right\}.
$$
Rather than giving the general definition of 2-nondegeneracy, we explain what this condition means in the case at hand.

Let $M$ be a 5-dimensional CR-hypersurface uniformly Levi degenerate of rank 1. Fix $p_0\in M$. Locally near $p_0$ the CR-structure is given by 1-forms $\mu$, $\eta^{\alpha}$, $\alpha=1,2$, where $\mu$ is\linebreak $i\RR$-valued and vanishes exactly on the complex tangent spaces $H_p(M)$, and $\eta^{\alpha}$ are $\CC$-valued and their restrictions to $H_p(M)$ at every point $p$ are $\CC$-linear and constitute a basis of $H_p^*(M)$. The integrability condition for the CR-structure is then equivalent to the Frobenius condition, which states that $d\mu$, $d\eta^{\alpha}$ belong to the differential ideal generated by $\mu$, $\eta^{\beta}$. Since $\mu$ is $i\RR$-valued, this implies
\begin{equation}
d\mu\equiv h_{\alpha\overline{\beta}}\eta^{\alpha}\wedge\eta^{\overline{\beta}}\quad(\mod \mu)\label{integr0}
\end{equation}
for some functions $h_{\alpha\overline{\beta}}$ satisfying $h_{\alpha\overline{\beta}}=h_{\overline{\beta}\alpha}$, where we use the summation convention for subscripts and superscripts (here and everywhere below the conjugation of indices denotes the conjugation of the corresponding forms, e.g. $\eta^{\bar\beta}:=\overline{\eta^{\beta}}$). Since $M$ is uniformly Levi degenerate of rank 1, one can choose $\eta^{\alpha}$ near $p_0$ so that 
\begin{equation}
(h_{\alpha\overline{\beta}})\equiv\left(
\begin{array}{cc}
\pm 1 & 0\\
0 & 0
\end{array}
\right).\label{hform}
\end{equation}
Further, the integrability condition yields
$$
d\eta^1\equiv \eta^2\wedge\sigma\quad(\mod \mu,\eta^1)
$$
for some complex-valued 1-form $\sigma$. Now, assuming that (\ref{hform}) holds, we say that $M$ is\linebreak 2-nondegenerate at $p_0$ if the coefficient at $\eta^{\bar 1}$ in the expansion of $\sigma$ with respect to $\mu$, $\eta^{\alpha}$, $\eta^{\bar\alpha}$ does not vanish at $p_0$. Clearly, with (\ref{hform}) satisfied, this condition is independent of the choice of $\mu$, $\eta^{\alpha}$. Finally, we say that $M$ is 2-nondegenerate if $M$ is 2-nondegenerate at every point. As shown in \cite{E} (see Proposition 1.16 and p. 51 therein), this definition of 2-nondegeneracy is equivalent to the standard one.    

Define ${\mathfrak C}_{2,1}$ to be the class of connected 5-dimensional CR-hypersurfaces that are 2-non\-dege\-nerate and uniformly Levi degenerate of rank 1. This class is quite large. Indeed, as explained in Section \ref{curvature}, the tube over the graph of a generic solution to the homogeneous Monge-Amp\`ere equation on $\RR^2$ lies in ${\mathfrak C}_{2,1}$. More examples are given by everywhere characteristic hypersurfaces for homogeneous differential operators and by tubes over homogeneous algebraic varieties (see Section 6 in \cite{E}). The main result of this paper is reduction of CR-structures in this class to absolute parallelisms (see Theorem \ref{main}). Such reduction was attempted earlier in \cite{E}, but the construction presented there has turned out to be only applicable to a more restricted class of ${\bf G}$-structures (see the correction to \cite{E}). Apart from article \cite{E}, we are not aware of any reduction results in the Levi degenerate case. We stress that our construction is fundamentally different from the one proposed in \cite{E}.      

We start by choosing a model of which any manifold in ${\mathfrak C}_{2,1}$ will be locally regarded as a deformation (see Section \ref{model}). For any CR-manifold $M$ let $\Aut(M)$ be the group of all CR-automorphisms of $M$, i.e.~CR-isomorphisms of $M$ onto itself. For $M$ to be chosen as a model, it needs to be homogeneous under $\Aut(M)$ and also locally the \lq\lq most symmetric\rq\rq. The latter condition means that at every point $p\in M$ the Lie algebra $\aut(M,p)$ of germs of infinitesimal CR-automorphisms of $M$ at $p$ must have the largest possible dimension, where an infinitesimal CR-automorphism is a smooth vector field whose flow consists of CR-maps. Consider the tube hypersurface over the future light cone in $\RR^3$:
\begin{equation}
M_0:=\bigl\{(z_1,z_2,z_3)\in\CC^3:(\Re z_1)^2+(\Re z_2)^2-(\Re z_3)^2=0,\,\Re z_3>0\bigr\}.\label{light}
\end{equation}
This hypersurface has been extensively studied (see, e.g. \cite{FK1}, \cite{FK2}, \cite{KZ}, \cite{Me}). In particular, $M_0$ is locally the most symmetric manifold among all locally homogeneous manifolds in ${\mathfrak C}_{2,1}$. More precisely, the algebra $\aut(M_0,p)$ at any point $p\in M_0$ is isomorphic to ${\mathfrak{so}}(3,2)$ (see \cite{KZ}) and has the largest dimension among all algebras of germs of infinitesimal CR-automorphisms of 5-dimen\-sion\-al locally homogeneous 2-nondegenerate CR-hypersurfaces (see Theorems I and II in \cite{FK2}). Furthermore, there is a CR-embedding of $M_0$ as an open dense subset in a homogeneous hypersurface $\Gamma$ in a complex projective quadric such that every local CR-automorphism of $M_0$ (i.e.~a CR-isomorphism between a pair of domains in $M_0$) continues to a global CR-automorphism of $\Gamma$ (see \cite{KZ}). These circumstances point to $\Gamma$ as a correct choice of a homogeneous model. We note that $\pi_1(\Gamma)\simeq\ZZ_2$ (see p. 69 in \cite{FK1}), and therefore another reasonable choice of a model would be the double cover $\hat\Gamma$ of $\Gamma$ (see Section \ref{applications} for details).

Next, in Section \ref{construction}, for every $M\in{\mathfrak C}_{2,1}$ we construct a fiber bundle ${\mathcal P}_M$ and an absolute parallelism $\omega_M$ on ${\mathcal P}_M$ with the required properties. Here ${\mathcal P}_M$ is a principal $H$-bundle, where $H$ is the isotropy subgroup of a point in $\Gamma$ under the action of the automorphism group\linebreak $G\simeq\SO(3,2)^{\circ}$ of $\Gamma$. In our construction we are guided by the Maurer-Cartan equations for $G$. Indeed, for $M=\Gamma$ the construction leads to the bundle $G\ra H\hspace{-0.1cm}\setminus\hspace{-0.1cm}G\simeq\Gamma$ and the right-invariant Maurer-Cartan form $\omega_G^{\hbox{\tiny MC}}$ on $G$, where $H\hspace{-0.1cm}\setminus\hspace{-0.1cm}G$ is the right coset space and $H$ acts on $G$ by left multiplication (generally, in our exposition principal bundles are realized as spaces with left actions of the structure groups). In Section \ref{construction} we also find a detailed expansion of the curvature $\Omega_M:=d\omega_M-\omega_M\wedge\omega_M$ of $\omega_M$. As shown in Section \ref{curvature} (see Theorem \ref{carconnec}), certain terms in this expansion (which we call the leading terms) are exactly the obstructions for $\omega_M$ to be a Cartan connection, i.e.~to change equivariantly under the action of $H$ on ${\mathcal P}_M$.     

Finally, in Section \ref{applications} we give applications of Theorem \ref{main}. For example, this theorem yields a new short proof of the extendability of local CR-automorphisms of $\Gamma$ and $\hat\Gamma$ to global CR-automorphisms mentioned above (see Corollary \ref{cor2}). Further, in Corollary \ref{cor4} we obtain the 2-jet determination property for germs of CR-isomorphisms between germs of manifolds in ${\mathfrak C}_{2,1}$. Next, in Corollary \ref{cor3} we show that for every $M\in{\mathfrak C}_{2,1}$ the group $\Aut(M)$ admits the structure of a Lie transformation group of $M$ of dimension at most 10, with dimension 10 occurring only for $\Gamma$ and $\hat\Gamma$. Finally, we settle Beloshapka's conjecture (see Conjecture \ref{belo}) for 5-dimensional manifolds in Corollary \ref{cor8}.

\section{The homogeneous model}\label{model}
\setcounter{equation}{0}

For $x=(x_1,\dots,x_5)\in\RR^5$ set
$$
(x,x):=x_1^2+x_2^2+x_3^2-x_4^2-x_5^2, 
$$
and let $\SO(3,2)$ be the group of all real $5\times 5$-matrices $C$ with $\det C=1$ satisfying $(Cx,Cx)\equiv (x,x)$, i.e.~$C^t{\mathbf J}C={\mathbf J}$ for
$$
{\mathbf J}:=\left(
\begin{array}{ll}
{\mathbf I}_3 & 0\\
0 & -{\mathbf I}_2
\end{array}
\right),
$$
where ${\mathbf I}_k$ is the identity $k\times k$-matrix. Consider the symmetric and Hermitian extensions of the form $(\,\,,\,)$ to $\CC^5$. Denote the symmetric extension by the same symbol $(\,\,,\,)$ and the Hermitian extension by $\langle\,\,,\,\rangle$. 

For $z=(z_1:\dots:z_5)\in\CC\PP^4$ we now consider the projective quadric
$$
{\mathcal Q}:=\{z\in\CC\PP^4: (z,z)=0\}
$$
and the open subset $\Omega\subset{\mathcal Q}$ defined as
\begin{equation}
\Omega:=\{z\in{\mathcal Q}:\langle z,z\rangle<0\}.\label{domainsomega}
\end{equation}
The set $\Omega$ has two connected components, one containing the point $p_{+}:=(0:0:0:1:i)$, the other containing the point $p_{-}:=(0:0:0:1:-i)$, and we denote these components by $\Omega_{\pm}$, respectively. For the following facts concerning $\Omega_{\pm}$ we refer the reader to pp. 285--289 in \cite{Sa}. First of all, we have
$$
\begin{array}{l}
\Omega_{+}=\left\{z\in\Omega: \Re z_4\Im z_5-\Re z_5\Im z_4>0\right\},\\
\vspace{-0.1cm}\\
\Omega_{-}=\left\{z\in\Omega:\Re z_4\Im z_5-\Re z_5\Im z_4<0\right\},
\end{array}
$$
and each of $\Omega_{\pm}$ is a realization of the symmetric classical domain of type $({\rm IV}_3)$, or, equivalently, of type $({\rm II}_2)$. Next, the group $\SO(3,2)$ acts on $\Omega$ effectively, and all holomorphic automorphisms of each of $\Omega_{\pm}$ arise from the action of the group ${\mathcal G}:=\SO(3,2)^{\circ}$, the connected identity component of $\SO(3,2)$. Observe also that the element
$$
g_0:=\left(\begin{array}{lllll}
0 & 1 & 0 & 0 & 0\\
1 & 0 & 0 & 0 & 0\\
0 & 0 & 1 & 0 & 0\\
0 & 0 & 0 & 0 & 1\\
0 & 0 & 0 & 1 & 0\\
\end{array}
\right)
$$
lies in the other connected component of $\SO(3,2)$ and interchanges the points $p_{\pm}$ hence the domains $\Omega_{\pm}$.

Further, it is not hard to see that the action of ${\mathcal G}$ on $\partial\Omega_{+}\cup\partial\Omega_{-}\subset{\mathcal Q}$ has two real hypersurface orbits
$$
\begin{array}{l}
\Gamma_{+}:=\left\{z\in{\mathcal Q}: (\Re z,\Re z)=(\Im z,\Im z)=(\Re z,\Im z)=0,\right.\\
\vspace{-0.3cm}\\
\hspace{1.05cm}\left.\Re z_4\Im z_5-\Re z_5\Im z_4>0\right\}\subset\partial\Omega_{+} \,\, \hbox{and}\\
\vspace{-0.1cm}\\
\Gamma_{-}:=\left\{z\in{\mathcal Q}: (\Re z,\Re z)=(\Im z,\Im z)=(\Re z,\Im z)=0,\right.\\
\vspace{-0.3cm}\\
\hspace{1.05cm}\left.\Re z_4\Im z_5-\Re z_5\Im z_4<0\right\}\subset\partial\Omega_{-}.
\end{array}
$$
Notice that for $q_{\pm}:=(\pm i:1: 0:1:\pm i)$ we have $\Gamma_{\pm}={\mathcal G}.q_{\pm}$. Clearly, the element $g_0$ of $\SO(3,2)$ interchanges the points $q_{\pm}$ and therefore the hypersurfaces $\Gamma_{\pm}$. Writing $\Omega_{\pm}$ in tube form (see \cite{FK1} and p. 289 in \cite{Sa}), one observes that the hypersurface $M_0$ introduced in (\ref{light}) is CR-equivalent to an open dense subset  of each of $\Gamma_{\pm}$, and since the hypersurfaces $\Gamma_{\pm}$ are homogeneous, it follows that they belong to the class ${\mathfrak C}_{2,1}$. By Theorems 4.5, 4.7 in \cite{KZ}, in this realization of $M_0$ its every local CR-automorphism extends to a global CR-automorphism of $\Gamma_{\pm}$ (in fact, to a holomorphic automorphism of ${\mathcal Q}$ induced by an element of $\SO(3,2)^{\circ}$). We note that although in \cite{KZ} only real-analytic maps were considered, every smooth local CR-automorphism of $M_0$ is in fact real-analytic, which is a consequence, for instance, of Theorem 7.1 in \cite{BJT}. 

Let $\Phi:(z_1:\dots:z_5)\mapsto(z_1^*:\dots:z_5^*)$ be the automorphism of $\CC\PP^4$ given by
$$
\begin{array}{lll}
\displaystyle z_1^*=\frac{1}{2}(z_1+iz_2-iz_4-z_5), & \displaystyle z_2^*=\frac{1}{2}(z_1-iz_2+iz_4-z_5), & z_3^*=z_3,\\
\vspace{-0.1cm}\\
\displaystyle z_4^*=\frac{1}{2}(z_1+iz_2+iz_4+z_5), & \displaystyle z_5^*=\frac{1}{2}(z_1-iz_2-iz_4+z_5).
\end{array}
$$
When viewed as a transformation of $\CC^5$, it takes $(\,\,,\,)$ and $\langle\,\,,\,\rangle$ into the bilinear and Hermitian forms defined, respectively, by the following matrices:
$$
{\mathbf S}:=\left(
\begin{array}{lllll}
0 & 0 & 0&0& 1\\
0 & 0 & 0&1& 0\\
0 & 0 & 1&0& 0\\
0 & 1 & 0&0& 0\\
1 & 0 & 0&0& 0\\
\end{array}
\right),\quad
{\mathbf T}:=\left(
\begin{array}{lllll}
0 & 0 & 0&1& 0\\
0 & 0 & 0&0& 1\\
0 & 0 & 1&0& 0\\
1 & 0 & 0&0& 0\\
0 & 1 & 0&0& 0\\
\end{array}
\right).
$$
Let $G$ be the connected identity component of the group of all complex $5\times 5$-matrices $C$ with $\det C=1$ satisfying
$$
C^t{\mathbf S}C = {\mathbf S}, \quad C^t{\mathbf T}\bar C ={\mathbf T}
$$
(clearly, $G$ is isomorphic to ${\mathcal G}$). The Lie algebra ${\mathfrak g}$ of $G$ (which is isomorphic to ${\mathfrak{so}}(3,2)$) consists of all matrices of the form
\begin{equation}
\left(\begin{array}{rrrrr}
\alpha & \beta & \gamma & \delta & 0\\
\bar \beta & \bar \alpha & \bar \gamma & 0 & -\delta\\
\sigma & \bar \sigma & 0 & -\bar\gamma & -\gamma\\
\rho & 0 & -\bar\sigma & -\bar\alpha & -\beta\\
0 & -\rho & -\sigma & -\bar\beta & -\alpha
\end{array}
\right),\label{liealgebra}
\end{equation}
with $\alpha,\beta,\gamma,\sigma\in\CC$, $\delta,\rho\in i\RR$. Set $\Gamma:=\Phi(\Gamma_{+})$, $q:=\Phi(q_{+})=(0:0:0:1:0)$. Then $\Gamma=G.q$, and we denote by $H\subset G$ the isotropy subgroup of $q$. We have $H=H^1\ltimes H^2$, where $H^1$ and $H^2$ are, respectively, the following subgroups of $H$:

\begin{equation}
\hspace{-1cm}\makebox[250pt]{$\left(
\begin{array}{ccccc}
A & 0 & 0&0& 0\\
0 & \bar A & 0&0& 0\\
0 & 0 & 1&0& 0\\
0 & 0 & 0&\bar A^{-1}& 0\\
0 & 0 & 0&0&A^{-1}\\
\end{array}
\right),\quad
\left(
\begin{array}{ccccc}
1 & 0 & 0&0& 0\\
0 & 1 & 0&0& 0\\
B & \bar B & 1&0& 0\\
\Lambda-|B|^2/2 & -\bar B^2/2 & -\bar B&1& 0\\
-B^2/2 & -\Lambda-|B|^2/2 & - B&0&1\\
\end{array}
\right),$}\label{subgroups}
\end{equation}
with $A\in\CC^*$, $B\in\CC$, $\Lambda\in i\RR$. 

Define a right action of $G$ on $\Gamma$ by
\begin{equation}
G\times\Gamma\ra\Gamma,\quad (g,p)\mapsto g^{-1}p,\label{action}
\end{equation}
and identify $\Gamma$ with the right coset space $H\hspace{-0.1cm}\setminus\hspace{-0.1cm}G$ by means of this action. Consider the principal $H$-bundle $G\ra H\hspace{-0.1cm}\setminus\hspace{-0.1cm}G\simeq\Gamma$ (with $H$ acting on $G$ by left multiplication) and the right-invariant Maurer-Cartan form $\omega_G^{\hbox{\tiny MC}}$ on $G$. In line with (\ref{liealgebra}), we now write $\omega_G^{\hbox{\tiny MC}}$ using scalar-valued 1-forms as follows:
\begin{equation}
\omega_G^{\hbox{\tiny MC}}=\left(
\begin{array}{ccccc}
\phi^2 & \theta^2 & \theta^1&\theta& 0\\
\theta^{\bar 2} & \phi^{\bar 2} &  \theta^{\bar 1}&0& -\theta\\
\phi^{\bar 1} & \phi^1 & 0&-\theta^{\bar 1}& -\theta^1\\
\psi & 0 & -\phi^1&- \phi^{\bar 2}& -\theta^2\\
0 & -\psi & -\phi^{\bar 1}&-\theta^{\bar 2}& -\phi^2\\
\end{array}\right),\label{mc}
\end{equation}
where $\theta^1$, $\theta^2$, $\phi^1$, $\phi^2$ are $\CC$-valued and $\theta,\phi$ take values in $i\RR$. The Maurer-Cartan equation $d\omega_G^{\hbox{\tiny MC}}=\omega_G^{\hbox{\tiny MC}}\wedge\omega_G^{\hbox{\tiny MC}}$ is then equivalent to the following set of identities:
\begin{equation}
\begin{array}{ll}
\hbox{(a)} & d\theta = -\theta^1\wedge \theta^{\bar 1} - \theta\wedge (\phi^2 +\phi^{\bar 2}),\\
\vspace{-0.1cm}\\
\hbox{(b)} & d\theta^1 = \theta^2 \wedge \theta^{\bar 1} -\theta^1\wedge \phi^2 - \theta\wedge \phi^1,\\
\vspace{-0.1cm}\\
\hbox{(c)} & d\theta^2 = -\theta^2\wedge (\phi^2 -\phi^{\bar 2}) +\theta^1 \wedge \phi^1,\\
\vspace{-0.1cm}\\
\hbox{(d)} & d\phi^1=-\theta^2\wedge\phi^{\bar 1}+\theta^1\wedge\psi+\phi^1\wedge\phi^{\bar 2},\\
\vspace{-0.1cm}\\
\hbox{(e)} & d\phi^2 =\theta^2\wedge\theta^{\bar 2}+\theta^1\wedge\phi^{\bar 1}+\theta\wedge \psi,\\
\vspace{-0.1cm}\\
\hbox{(f)} & d\psi=-\phi^1\wedge\phi^{\bar 1}+\psi\wedge(\phi^2+\phi^{\bar 2}).
\end{array}\label{structure1}
\end{equation}

In the next section to every $M\in{\mathfrak C}_{2,1}$ we associate a principal $H$-bundle ${\mathcal P}_M\ra M$ (with a left action of $H$ on ${\mathcal P}_M$) and a ${\mathfrak g}$-valued absolute parallelism $\omega_M$ on ${\mathcal P}_M$ such that the bundles ${\mathcal P}_{\Gamma}\ra\Gamma$ and $G\ra H\hspace{-0.1cm}\setminus\hspace{-0.1cm}G\simeq\Gamma$ are isomorphic and under this isomorphism the parallelisms $\omega_{\Gamma}$ and $\omega_G^{\hbox{\tiny MC}}$ are identified. The CR-manifold $\Gamma$ will be regarded as a \lq\lq flat\rq\rq\, model for the class ${\mathfrak C}_{2,1}$, and structure equations (\ref{structure1}) will guide us through our construction.

\section{The general case}\label{construction}
\setcounter{equation}{0}

Fix $M\in{\mathfrak C}_{2,1}$. As we mentioned in the introduction, locally on $M$ the CR-structure is given by 1-forms $\mu$, $\eta^{\alpha}$, $\alpha=1,2$, where $\mu$ is $i\RR$-valued and vanishes exactly on the complex tangent spaces $H_p(M)$, and $\eta^{\alpha}$ are $\CC$-valued and their restrictions to $H_p(M)$ at every point $p$ are $\CC$-linear and form a basis of $H_p^*(M)$. The integrability condition for the CR-structure then implies that identity (\ref{integr0}) holds for some functions $h_{\alpha\overline{\beta}}$ satisfying $h_{\alpha\overline{\beta}}=h_{\overline{\beta}\alpha}$.

For $p\in M$ define $E_p$ as the collection of all $i\RR$-valued covectors $\theta$ on $T_p(M)$ such that $H_p(M)=\{X\in T_p(M): \theta(X)=0\}$. Clearly, all elements in $E_p$ are real non-zero multiples of each other. Let $E$ be the bundle over $M$ with fibers $E_p$. Define $\omega$ to be the tautological 1-form on $E$, that is, for $\theta\in E$ and $X\in T_{\theta}\left(E\right)$ set
$$
\omega(\theta)(X):=\theta(d\pi_E(\theta)(X)),
$$
where $\pi_E: E\ra M$ is the projection. Since the Levi form of $M$ has rank 1 everywhere, identity (\ref{integr0}) implies that for every $\theta\in E$ there exist a real-valued covector $\phi$ and a complex-valued covector $\theta^1$ on $T_{\theta}(E)$ such that: (i) $\theta^1$ is the pull-back of a complex-valued covector on $T_{\pi_{{}_E}(\theta)}(M)$ complex-linear on $H_{\pi_{{}_E}(\theta)}(M)$, and (ii) the following identity holds:
\begin{equation}
d\omega(\theta)=\pm\theta^1\wedge\theta^{\overline{1}}-\omega(\theta)\wedge\phi.\label{integr}
\end{equation}

For every $p\in M$ the fiber $E_p$ has exactly two connected components, and the signs in the right-hand side of (\ref{integr}) coincide for all covectors $\theta$ lying in the same connected component of $E_p$ and are opposite for any two covectors lying in different connected components irrespectively of the choice of $\theta^1$, $\phi$. We then define a bundle ${\mathcal P}^1$ over $M$ as follows: for every $p\in M$ the fiber ${\mathcal P}^1_p$ over $p$ is connected and consists of all elements $\theta\in E_p$ for which the minus sign occurs in the right-hand side of (\ref{integr}); we also set $\pi^1:=\pi_E\bigr|_{{\mathcal P}^1}$. 

Next, the most general transformation of $(\omega(\theta),\theta^1,\theta^{\overline{1}},\phi)$ preserving the equation 
\begin{equation}
d\omega(\theta)=-\theta^1\wedge\theta^{\overline{1}}-\omega(\theta)\wedge\phi\label{integrminus}
\end{equation}
and the covector $\omega(\theta)$ is  given by the matrix (acting on the left)
\begin{equation}
\left(
\begin{array}{rccc}
1 & 0 & 0 & 0\\
\vspace{-0.3cm}\\
\overline{b} & a & 0 & 0\\
\vspace{-0.3cm}\\
-b & 0 & \overline{a} & 0\\
\vspace{-0.3cm}\\
\lambda & -ab & -\overline{a}\overline{b} & 1
\end{array}
\right),\label{g1structure}
\end{equation}
where $a,b\in\CC$, $|a|=1$, $\lambda\in i\RR$. Let $H_1$ be the group of matrices of the form (\ref{g1structure}). Observe that $H_1$ is isomorphic to the subgroup $H_1^1\ltimes H^2$ of $H$, where $H_1^1$ is the subgroup of $H^1$ given by the condition $|A|=1$ (see (\ref{subgroups})). Our goal is to reduce the $H_1$-structure on ${\mathcal P}^1$ to an absolute parallelism.

We now introduce a principal $H_1$-bundle ${\mathcal P}^2$ over ${\mathcal P}^1$ as follows: for $\theta\in {\mathcal P}^1$ let the fiber ${\mathcal P}_{\theta}^2$ over $\theta$ be the collection of all 4-tuples of covectors $(\omega(\theta),\theta^1,\theta^{\overline{1}},\phi)$ on $T_{\theta}({\mathcal P}^1)$, where $\theta^1$ and $\phi$ are chosen as described above. Let $\pi^2: {\mathcal P}^2\ra {\mathcal P}^1$ be the projection. It is easy to see that ${\mathcal P}^2$ is a principal $H$-bundle if considered as a fiber bundle over $M$ with projection $\pi:=\pi^1\circ\pi^2$. Indeed, using the subgroups $H^1$, $H^2$ introduced in (\ref{subgroups}) we define a (left) action of $H$ on ${\mathcal P}^2$ by the formulas

\begin{equation}
\hspace{-1cm}\makebox[250pt]{$\begin{array}{l}
(\theta,\theta^1,\theta^{\overline{1}},\phi)\mapsto (|A|^2\theta,A\theta^1,\overline{A}\theta^{\overline{1}},\phi),\\
\vspace{-0.1cm}\\
(\theta,\theta^1,\theta^{\overline{1}},\phi)\mapsto (\theta,\theta^1+\overline{B}\omega(\theta),\theta^{\overline{1}}-B\omega(\theta),\phi-B\theta^1-\overline{B}\theta^{\overline{1}}-2\Lambda\omega(\theta)),
\end{array}$}\label{actionbyh}
\end{equation}
where we identified any two covectors $\theta^1$ and $\theta^{1'}$ at points $\theta$ and $\theta'$ of a fiber ${\mathcal P}_p^1$ that are obtained by pulling back the same covector at $p\in M$.

We now define two tautological 1-forms on ${\mathcal P}^2$ as
$$
\begin{array}{l}
\omega^1({\bf\Theta})(X):=\theta^1(d\pi^2({\bf\Theta})(X)),\\
\vspace{-0.3cm}\\
\varphi({\bf\Theta})(X):=\phi(d\pi^2({\bf\Theta})(X)),
\end{array}
$$
where ${\bf\Theta}=(\omega(\theta),\theta^1,\theta^{\overline{1}},\phi)$ is a point in ${\mathcal P}_{\theta}^2$ and $X\in T_{{\bf\Theta}}({\mathcal P}^2)$. It is clear from (\ref{integrminus}) that these forms satisfy
\begin{equation}
d\omega=-\omega^1\wedge\omega^{\overline{1}}-\omega\wedge\varphi,\label{integrplushigh}
\end{equation}
where we denote the pull-back of $\omega$ from ${\mathcal P}^1$ to ${\mathcal P}^2$ by the same symbol (cf. equation (a) in (\ref{structure1})). Further, computing $d\omega^1$ in local coordinates on ${\mathcal P}^2$ and using the integrability of the CR-structure of $M$ we obtain
\begin{equation}
d\omega^1=\theta^{2}\wedge\xi-\omega^1\wedge\varphi^2-\omega\wedge\varphi^{1}\label{integr1}
\end{equation}
for some complex-valued 1-forms $\theta^2,\xi,\varphi^1,\varphi^2$ (cf. equation (b) in (\ref{structure1})). Here for any\linebreak ${\bf\Theta}=(\omega(\theta),\theta^1,\theta^{\overline{1}},\phi)$ the covector $\theta^2({\bf\Theta})$ is the pull-back of a complex-valued covector $\theta_0^2$ at $p:=\pi({\bf\Theta})$ such that $\theta_0^2$ is complex-linear on $H_p(M)$ and the restrictions of $\theta_0^1$ and $\theta_0^2$ to $H_p(M)$ form a basis of $H_p^*(M)$, where $\theta_0^1$ is the covector at $p$ that pulls back to $\theta^1$. 

In what follows we study consequences of identities (\ref{integrplushigh}) and (\ref{integr1}). Our calculations will be entirely local, and we will impose conditions that will determine the forms $\theta^2,\varphi^1,\varphi^2$ (as well as another $i\RR$-valued 1-form $\psi$ introduced below) uniquely. This will allow us to patch the locally defined forms $\theta^2,\varphi^1,\varphi^2$, $\psi$ into globally defined 1-forms on ${\mathcal P}^2$. Together with $\omega$, $\omega^1$, $\varphi$ these globally defined forms will be used to construct an absolute parallelism on ${\mathcal P}^2$.

Exterior differentiation of (\ref{integrplushigh}) and substitution of (\ref{integrplushigh}), (\ref{integr1}) for $d\omega$, $d\omega^1$, respectively, yield
\begin{equation}
\begin{array}{l}
(\varphi-\varphi^2-\varphi^{\bar 2})\wedge\omega^1\wedge\omega^{\bar 1}+{\bar\xi}\wedge\omega^1\wedge\theta^{\bar 2}+\xi\wedge\theta^2\wedge\omega^{\bar 1}+\\
\vspace{-0.1cm}\\
\hspace{7cm}(d\varphi-\omega^1\wedge\varphi^{\bar 1}-\omega^{\bar 1}\wedge\varphi^1)\wedge\omega=0.
\end{array}\label{eq1}
\end{equation}
It then follows from Cartan's lemma that $\varphi-\varphi^2-\varphi^{\bar 2}=P\omega^1+\overline{P}\omega^{\bar 1}+Q\theta^2+\overline{Q}\theta^{\bar 2}+R\omega$ for some smooth functions $P,Q,R$, where $R$ is $i\RR$-valued. Setting
$$
\tilde\varphi^2:=\varphi^2+P\omega^1+Q\theta^2+\frac{1}{2}R\omega,
$$
we see that the form $\varphi^2$ can be chosen to satisfy identity (\ref{integr1}) with some $\tilde\xi$, $\tilde\varphi^1$ in place of $\xi$, $\varphi^1$ as well as the condition
\begin{equation}
\Re\varphi^2=\frac{\varphi}{2},\label{cond1}
\end{equation}
and from now on we assume that (\ref{cond1}) holds (cf. equations (a) and (b) in (\ref{structure1})). Note that (\ref{cond1}) is analogous to normalization condition (4.21) in \cite{CM}. 

With conditions (\ref{integr1}), (\ref{cond1}) satisfied, identity (\ref{eq1}) and Cartan's lemma imply\linebreak $\xi=U\theta^2+V\omega^{\bar 1}+W\omega$ for some functions $U,V,W$. Setting $\tilde\xi:=\xi-U\theta^2-W\omega$, we therefore can assume that the form $\xi$ is a multiple of $\omega^{\bar 1}$ and satisfies (\ref{integr1}) with some $\tilde\varphi^1$ in place of $\varphi^1$. For $\xi=V\omega^{\bar 1}$ the 2-nondegeneracy of $M$ yields that $V$ is a nowhere vanishing function, thus by scaling $\theta^2$ we suppose from now on that $\xi=\omega^{\bar 1}$. Hence identity (\ref{integr1}) turns into the identity
\begin{equation}
d\omega^1=\theta^{2}\wedge\omega^{\bar 1}-\omega^1\wedge\varphi^2-\omega\wedge\varphi^{1}\label{integr11}
\end{equation}
(cf. equation (b) in (\ref{structure1})). Furthermore, (\ref{eq1}) now implies
\begin{equation}
d\varphi=\omega^1\wedge\varphi^{\bar 1}+\omega^{\bar 1}\wedge\varphi^1+2\omega\wedge\psi\label{integr2}
\end{equation}
for some $i\RR$-valued form $\psi$ (cf. equation (e) in (\ref{structure1})).    

Next, it is not hard to see that the forms $\theta^2$, $\varphi^1$, $\varphi^2$, $\psi$ satisfying (\ref{cond1}), (\ref{integr11}), (\ref{integr2}) are defined up to the following transformations:
\begin{equation}
\begin{array}{l}
\theta^2=\tilde\theta^2+c\omega^1+f\omega,\\
\vspace{-0.3cm}\\
\varphi^2=\tilde\varphi^2-\overline{c}\omega^1+c\omega^{\bar 1}+g\omega,\\
\vspace{-0.3cm}\\
\varphi^1=\tilde\varphi^1+g\omega^1+f\omega^{\bar 1}+r\omega,\\
\vspace{-0.3cm}\\
\displaystyle\psi=\tilde\psi-\frac{\overline{r}}{2}\omega^1+\frac{r}{2}\omega^{\bar 1}+s\omega,
\end{array}\label{transform}
\end{equation}
for some functions $c$, $f$, $g$, $r$, $s$, where $g$ and $s$ are $\RR$-valued. Indeed, the first equation in (\ref{transform}) is the most general change of $\theta^2$ preserving the form of (\ref{integr11}). Next, subtracting (\ref{integr11}) from the same identity
with $\tilde\theta^2$, $\tilde\varphi^1$, $\tilde\varphi^2$ in place of $\theta^2$, $\varphi^1$, $\varphi^2$, we obtain
$$
-(c\omega^1+ f\omega)\wedge\omega^{\bar 1}-
\omega^1\wedge(\tilde\varphi^2-\varphi^2)-
\omega\wedge (\tilde\varphi^1 - \varphi^1)=0,
$$
which implies, together with (\ref{cond1}) and Cartan's lemma, the second and third equations in (\ref{transform}). Finally, subtracting (\ref{integr2})  from the same identity
with  $\tilde\varphi^1$, $\tilde\psi$ in place of $\varphi^1$, $\psi$, we obtain
$$
\omega^1\wedge(-g\omega^{\bar 1}-
\bar f\omega^{1}+\bar r\omega)+\omega^{\bar 1}\wedge (-g\omega^1-f\omega^{\bar 1}-r\omega)+
2\omega\wedge(\tilde\psi - \psi)=0,
$$
which leads to the last equation in (\ref{transform}). 

As computation in local coordinates immediately shows, for any forms $\theta^2$, $\varphi^1$, $\varphi^2$, $\psi$ satisfying (\ref{cond1}), (\ref{integr11}), (\ref{integr2}) the values of $i\omega$, $\Re\omega^1$, $\Im\omega^1$, $\Re\theta^2$, $\Im\theta^2$, $\Re\varphi^1$, $\Im\varphi^1$, $\varphi$, $\Im\varphi^2$, $i\psi$ at any ${\bf\Theta}$ constitute a basis of $T_{{\bf\Theta}}^*({\mathcal P}^2)$. In what follows we utilize expansions of complex-valued forms on ${\mathcal P}^2$ with respect to  $\omega$, $\omega^1$, $\omega^{\bar 1}$, $\theta^2$, $\theta^{\bar 2}$, $\varphi^1$, $\varphi^{\bar 1}$, $\varphi^2$, $\varphi^{\bar 2}$, $\psi$. We will be particularly interested in coefficients at wedge products of $\omega$, $\omega^1$, $\omega^{\bar 1}$, $\theta^2$, $\theta^{\bar 2}$ and for a form $\Omega$ denote them by $\Omega_{\alpha\dots{\bar\beta}\dots\,0}$, where $\alpha,\beta=1,2$, with index 0 corresponding to $\omega$, index 1 to $\omega^1$, and index 2 to $\theta^2$.  

We will now impose further conditions on the forms $\theta^2$, $\varphi^1$, $\varphi^2$, $\psi$ in order to fix them uniquely. Our conditions are inspired by equations (c), (d), (e), (f) in (\ref{structure1}). First, let
\begin{equation}
\Theta^2:=d\theta^2+\theta^2\wedge(\varphi^2-\varphi^{\bar 2})-\omega^1\wedge\varphi^1\label{Sigma}
\end{equation}
(cf. equation (c) in (\ref{structure1})). By the integrability of the CR-structure on $M$ we have
\begin{equation}
d\theta^2=\theta^2\wedge\zeta^2+\omega^1\wedge\zeta^1+\omega\wedge\zeta\label{integr3}
\end{equation}
for some complex-valued 1-forms $\zeta^1$, $\zeta^2$, $\zeta$. Plugging (\ref{integr3}) into (\ref{Sigma}) we obtain
\begin{equation}
\Theta^2=\theta^2\wedge(\zeta^2+\varphi^2-\varphi^{\bar 2})+\omega^1\wedge(\zeta^1-\varphi^1)+\omega\wedge\zeta.\label{theta2}
\end{equation}
First of all, we need transformation formulas for $\zeta^1$, $\zeta^2$ when $\theta^2$ changes as in (\ref{transform}).

\begin{lemma}\label{lemmatransf}\sl For any forms  $\tilde\zeta^1$, $\tilde\zeta^2$ corresponding to $\tilde\theta^2$ one has
\begin{equation}
\begin{array}{ll}
\tilde\zeta^2\equiv\zeta^2-c\omega^{\bar 1}& (\hbox{\rm mod}\,\omega,\omega^1, \theta^2),\\
\vspace{-0.1cm}\\
\tilde\zeta^1\equiv \zeta^1+f\omega^{\bar 1}-c^2\omega^{\bar 1}+c\varphi^2+c\zeta^2+dc& (\hbox{\rm mod}\,\omega, \omega^1,\theta^2).
\end{array}\label{eqq1}
\end{equation}
\end{lemma}

\noindent{\bf Proof:} Using (\ref{integrplushigh}), (\ref{integr11}), (\ref{transform}), (\ref{integr3}) we calculate
$$
\begin{array}{l}
\tilde\theta^2\wedge\tilde\zeta^2+\omega^1\wedge\tilde\zeta^1+\omega\wedge\tilde\zeta=d\tilde\theta^2=d(\theta^2-c\omega^1-f\omega)=\\
\vspace{-0.1cm}\\
d\theta^2-dc\wedge\omega^1-cd\omega^1-df\wedge\omega-fd\omega=\theta^2\wedge\zeta^2+\omega^1\wedge\zeta^1+\omega\wedge\zeta-\\
\vspace{-0.1cm}\\
dc\wedge\omega^1-c(\theta^2\wedge\omega^{\bar 1}-\omega^1\wedge\varphi^2-\omega\wedge\varphi^1)-df\wedge\omega-f(-\omega^1\wedge\omega^{\bar 1}-\omega\wedge\varphi)=\\
\vspace{-0.1cm}\\
\theta^2\wedge(\zeta^2-c\omega^{\bar 1})+\omega^1\wedge(\zeta^1+dc+c\varphi^2+f\omega^{\bar 1})+\omega\wedge(\zeta+c\varphi^1+df+f\varphi)=\\
\vspace{-0.1cm}\\
(\tilde\theta^2+c\omega^1+f\omega)\wedge(\zeta^2-c\omega^{\bar 1})+\omega^1\wedge(\zeta^1+dc+c\varphi^2+f\omega^{\bar 1})+\\
\vspace{-0.1cm}\\
\omega\wedge(\zeta+c\varphi^1+df+f\varphi)=\\
\vspace{-0.1cm}\\
\tilde\theta^2\wedge(\zeta^2-c\omega^{\bar 1})+\omega^1\wedge(c\zeta^2-c^2\omega^{\bar 1}+\zeta^1+dc+c\varphi^2+f\omega^{\bar 1})+\\
\vspace{-0.1cm}\\
\omega\wedge(f\zeta^2-fc\omega^{\bar 1}+\zeta+c\varphi^1+df+f\varphi).
\end{array}
$$
Formulas (\ref{eqq1}) now follow from Cartan's lemma.\qed
\vspace{-0.3cm}\\

We will fix the functions $c$, $f$ by imposing conditions on  $\Theta^2_{2{\bar1}}$, $\Theta^2_{1{\bar 1}}$. Thus, we need to know how these coefficients are transformed under (\ref{transform}). First we consider $\Theta^2_{2{\bar1}}$. Exterior differentiation of  (\ref{integr11}) and substitution of  (\ref{integrplushigh}), (\ref{integr11}),  (\ref{integr3}) for $d\omega$, $d\omega^1$, $d\theta^2$, respectively, yield
\begin{equation}
\hspace{-1cm}\makebox[250pt]{$\begin{array}{l}
(-\zeta^2-\varphi^2+\varphi^{\bar 2})\wedge\theta^2\wedge\omega^{\bar 1}+((\zeta^1-\varphi^1)\wedge\omega^{\bar 1}+d\varphi^2-\theta^2\wedge\theta^{\bar 2})\wedge\omega^1+\\
\vspace{-0.1cm}\\
\hspace{7cm}(\zeta\wedge\omega^{\bar 1}+d\varphi^1+\theta^2\wedge\varphi^{\bar 1}-\varphi^1\wedge\varphi^{\bar 2})\wedge\omega=0.
\end{array}$}\label{eqq2}
\end{equation}
By Cartan's lemma, formula (\ref{eqq2}) implies 
\begin{equation}
\zeta^2+\varphi^2-\varphi^{\bar 2}\equiv 0\quad(\mod \omega,\omega^1,\omega^{\bar 1},\theta^2),\label{normal8}
\end{equation}
hence by (\ref{theta2}) we have
$$
\Theta^2 \equiv \Theta^2_{2\bar 1}\theta^2\wedge\omega^{\bar1}\quad (\mod \omega, \omega^1).
$$
Further, it follows from (\ref{transform}) and the first equation in (\ref{eqq1}) that 
$$
\tilde\zeta^2+\tilde\varphi^2-\tilde\varphi^{\bar 2}\equiv (\zeta^2+\varphi^2-\varphi^{\bar 2})-3c\omega^{\bar 1}\, (\hbox{\rm mod}\,\omega,\omega^1, \theta^2).
$$
Therefore, passing to the forms with tildas we see
\begin{equation}
\tilde\Theta^2_{2{\bar1}}=\Theta^2_{2{\bar1}}-3c,\label{choiceofc}
\end{equation}
which shows that $c$ is uniquely fixed by the requirement $\tilde\Theta^2_{2{\bar1}}=0$. Thus, from now on we assume that the condition
\begin{equation}
\Theta^2_{2{\bar1}}=0\label{curv1}
\end{equation}
is satisfied, or, equivalently, $\Theta^2 \equiv 0\,(\mod \omega, \omega^1).$
Hence $\zeta^2+\varphi^2-\varphi^{\bar 2}\equiv 0\,(\hbox{\rm mod}\,\omega,\omega^1,\theta^2)$ and in formula (\ref{transform}) one has $c=0$.

Now, by choosing suitable $\zeta^1$ and $\zeta$ in (\ref{integr3}),  
we can assume that
\begin{equation}
\zeta^2=-\varphi^2+\varphi^{\bar 2},\label{aux2}
\end{equation}
hence (\ref{eqq2}) together with Cartan's lemma implies
\begin{equation}
(\zeta^1-\varphi^1)\wedge\omega^{\bar 1}+d\varphi^2-\theta^2\wedge\theta^{\bar 2} \equiv 0\quad (\mod \omega,\omega^1).\label{dphi}
\end{equation}
By (\ref{cond1}), (\ref{integr2}), (\ref{aux2}) we also have
\begin{equation}
d\zeta^2=-d\varphi^2+d\varphi^{\bar2}= d\varphi - 2d\varphi^2
\equiv \omega^{\bar1}\wedge\varphi^1 - 2d\varphi^2\quad(\mod \omega,\omega^1).\label{aux3}
\end{equation}
Further, exterior differentiation of  (\ref{integr3}) and substitution of (\ref{integrplushigh}), (\ref{integr11}), (\ref{integr3}) for $d\omega$, $d\omega^1$, $d\theta^2$, respectively, lead to
\begin{equation}
\hspace{-0.15cm}\begin{array}{l}
-(\zeta^1\wedge\omega^{\bar 1}+d\zeta^2)\wedge\theta^2+(\zeta^1\wedge\zeta^2+\zeta^1\wedge\varphi^2+\zeta\wedge\omega^{\bar 1}-d\zeta^1)\wedge\omega^1+\\
\vspace{-0.1cm}\\
\hspace{8cm}(\zeta\wedge\zeta^2+\zeta^1\wedge\varphi^1+\zeta\wedge\varphi-d\zeta)\wedge\omega=0.
\end{array}\label{exp}
\end{equation}
In particular, Cartan's lemma implies
$$
\zeta^1\wedge\omega^{\bar 1}+d\zeta^2 \equiv 0\quad (\mod \omega,\omega^1, \theta^2),
$$
hence by (\ref{aux3})
$$
(\zeta^1-\varphi^1)\wedge\omega^{\bar1}-2d\varphi^2\equiv0\quad (\mod \omega, \omega^1,\theta^2).
$$
After eliminating $d\varphi^2$ by means of (\ref{dphi}) this identity leads to
$$
(\zeta^1-\varphi^1)\wedge\omega^{\bar1} \equiv 0\quad(\mod \omega,\omega^1, \theta^2),
$$
and therefore to
$$
\zeta^1-\varphi^1 \equiv 0\quad(\mod \omega,\omega^1, \omega^{\bar 1},\theta^2).
$$
Then 
$$
\Theta^2= \omega^1\wedge (\zeta^1 -\varphi^1) + \omega\wedge\zeta\equiv
\Theta^2_{1\bar1} \omega^1\wedge\omega^{\bar1}\quad(\mod\omega, \theta^2\wedge\omega^1).
$$
Further, it follows from (\ref{transform}) and the second equation in (\ref{eqq1}) for $c=0$ that 
$$
\tilde\zeta^1-\tilde\varphi^1\equiv (\zeta^1-\varphi^1)+2f\omega^{\bar 1}\, (\hbox{\rm mod}\,\omega,\omega^1, \theta^2).
$$
Therefore, passing to the forms with tildas we see
$$
\tilde\Theta^2_{1{\bar1}}=\Theta^2_{1{\bar1}}+2f,
$$
which shows that $f$ is uniquely fixed by the requirement $\tilde\Theta^2_{1{\bar1}}=0$. Thus, from now on we assume that the condition
\begin{equation}
\Theta^2_{1{\bar1}}=0\label{curv2}
\end{equation}
is satisfied. Hence $\zeta^1-\varphi^1\equiv0\,(\hbox{\rm mod}\,\omega,\omega^1, \theta^2)$ and in formula (\ref{transform}) one has $c=f=0$. Observe that $\theta^2$ has now been canonically fixed. 

Next, by adjusting $\zeta$ in (\ref{integr3}) we assume 
$$
\zeta^1\equiv \varphi^1\quad (\mod \theta^2),
$$
i.e.
$$
\Theta^2 = \Theta^2_{21}\theta^2\wedge\omega^1 + \omega\wedge\zeta,
$$
and
\begin{equation}
\zeta^1 = \varphi^1 -  \Theta^2_{21}\theta^2.\label{aux5}
\end{equation}
Further, from (\ref{exp}) using Cartan's lemma we obtain
$$
\zeta^1\wedge\zeta^2+\zeta^1\wedge\varphi^2+\zeta\wedge\omega^{\bar 1}-d\zeta^1\equiv 0\quad (\mod\omega,\omega^1, \theta^2),
$$
i.e.~by (\ref{aux2}), (\ref{aux5})
$$
(\varphi^1 -  \Theta^2_{21}\theta^2)\wedge\varphi^{\bar2}+\zeta\wedge\omega^{\bar 1}-d\zeta^1\equiv 0\quad  (\mod\omega,\omega^1, \theta^2).
$$
Hence using (\ref{aux5}) once again we see
\begin{equation}
\zeta\wedge \omega^{\bar1} =
(\Theta^2_{21}\theta^2  -\varphi^1)\wedge\varphi^{\bar2}
+ d\varphi^1  - d(\Theta^2_{21}\theta^2)
\equiv -\varphi^1\wedge\varphi^{\bar2}
+ d\varphi^1\quad (\mod\omega,\omega^1, \theta^2).\label{aux6}
\end{equation}
On the other hand, by (\ref{eqq2}), (\ref{aux2}) we have
\begin{equation}
\zeta\wedge\omega^{\bar 1}+d\varphi^1+\theta^2\wedge\varphi^{\bar 1}-\varphi^1\wedge\varphi^{\bar 2}\equiv 0\quad (\mod \omega,\omega^1),\label{aux77}
\end{equation}
and by adding (\ref{aux6}), (\ref{aux77})
we obtain
$$
\zeta\wedge\omega^{\bar 1}\equiv 0\quad (\mod \omega,\omega^1,\theta^2),
$$
which yields
$$
\zeta \equiv 0\quad (\mod \omega,\omega^1,\omega^{\bar 1},\theta^2).
$$
Thus, we have obtained the following expansion for the 2-form $\Theta^2$:
\begin{equation}
\Theta^2 = \Theta^2_{21}\theta^2\wedge\omega^1 
+\Theta^2_{20} \theta^2\wedge\omega 
+\Theta^2_{10}\omega^1\wedge\omega
+\Theta^2_{\bar 10} \omega^{\bar 1}\wedge\omega.\label{curvv1}
\end{equation}
Accordingly, in (\ref{integr3}) we choose
\begin{equation}
\zeta =  - \Theta^2_{20}\theta^2-\Theta^2_{10}\omega^1-\Theta^2_{\bar10}\omega^{\bar1}.\label{aux7}
\end{equation}

We will now fix the function $g$. For this purpose we use the following form:
\begin{equation}
\hspace{-0.3cm}\Phi^2:=d\varphi^2-\theta^2\wedge\theta^{\bar 2}-\omega^1\wedge\varphi^{\bar 1}-\omega\wedge\psi\label{Phi2}
\end{equation}
(cf. equation (e) in (\ref{structure1})). Observe that according to (\ref{cond1}), (\ref{integr2}) we have $\Re\Phi^2=0$. Formulas (\ref{dphi}), (\ref{aux5}) yield
$$
d\varphi^2 \equiv \theta^2\wedge\theta^{\bar2} +  \Theta^2_{21}\theta^2\wedge\omega^{\bar1}\quad (\mod\omega, \omega^1),\label{aux8}
$$
hence
$$
\Phi^2\equiv\Theta^2_{21}\theta^2\wedge\omega^{\bar 1}\quad(\mod\omega,\omega^1).
$$
Since $\Re\Phi^2=0$, we obtain
\begin{equation}
\Phi^2=\Theta^2_{21}\theta^2\wedge\omega^{\bar 1}+\Theta^{\bar 2}_{{\bar 2}{\bar 1}}\omega^1\wedge\theta^{\bar 2}+\Phi^2_{1\bar 1}\omega^1\wedge\omega^{\bar 1}+\Phi^2_0\wedge\omega,\label{aux9}
\end{equation}
where $\Phi^2_{1\bar 1}$ is a real-valued function and $\Phi^2_0$ is a real-valued 1-form. Now, by (\ref{transform}), (\ref{Phi2}) we have
$$
\tilde\Phi^2\equiv \Phi^2+2g\omega^1\wedge\omega^{\bar 1}\quad(\mod\omega).
$$
Since the expansion of $\Phi^2$ in (\ref{aux9}) does not involve $\varphi^1\wedge\psi$, $\varphi^{\bar 1}\wedge\psi$, $\omega^1\wedge\psi$, $\omega^{\bar 1}\wedge\psi$, $\omega^1\wedge\varphi^{\bar 1}$, $\omega^{\bar 1}\wedge\varphi^1$, $\varphi^1\wedge\varphi^{\bar 1}$, it follows that
$$
\tilde\Phi^2_{1\bar 1}=\Phi^2_{1\bar 1}+2g,
$$
which shows that $g$ is uniquely fixed by the requirement $\tilde\Phi^2_{1\bar 1}=0$ (recall that $\Phi^2_{1\bar 1}$ is real-valued). Thus, from now on we assume that the condition
\begin{equation}
\Phi^2_{1\bar 1}=0\label{curv3}
\end{equation}
is satisfied, hence in formula (\ref{transform}) one has $c=f=g=0$. Observe that $\varphi^2$ has now been canonically fixed.

We will now fix the function $r$. For this purpose we use the following form:
\begin{equation}
\Phi^1:=d\varphi^1+\theta^2\wedge\varphi^{\bar 1}-\omega^1\wedge\psi-\varphi^1\wedge\varphi^{\bar 2}\label{Phi1}
\end{equation}
(cf. equation (d) in (\ref{structure1})). Substitution of (\ref{aux2}), (\ref{aux5}), (\ref{aux7}) into (\ref{eqq2}) implies
$$
\begin{array}{l}
(-\Theta^2_{21}\theta^2\wedge\omega^{\bar 1}+d\varphi^2-\theta^2\wedge\theta^{\bar 2})\wedge\omega^1+\\
\vspace{-0.1cm}\\
\hspace{5cm}(- \Theta^2_{20}\theta^2 \wedge\omega^{\bar 1}-\Theta^2_{10}\omega^1\wedge\omega^{\bar 1}+d\varphi^1+\theta^2\wedge\varphi^{\bar 1}-\varphi^1\wedge\varphi^{\bar2})\wedge\omega=0
\end{array}
$$
and together with (\ref{Phi2}), (\ref{aux9}) leads to
\begin{equation}
d\varphi^1\equiv-\theta^2\wedge\varphi^{\bar 1}+\varphi^1\wedge\varphi^{\bar 2}+\Theta^2_{20}\theta^2 \wedge\omega^{\bar 1}+\Theta^2_{10}\omega^1\wedge\omega^{\bar 1}+(\Phi^2_0-\psi)\wedge\omega^1\quad(\mod\omega),\label{aux18}
\end{equation}
hence
\begin{equation}
\Phi^1=\Theta^2_{20}\theta^2 \wedge\omega^{\bar 1}+(\Phi^2_0-\Theta^2_{10}\omega^{\bar 1})\wedge\omega^1+\Phi^1_0\wedge\omega\label{aux10}
\end{equation}
for some 1-form $\Phi^1_0$. 

To proceed further, we need to know what terms are involved in the expansion of $\Phi^2_0$ in (\ref{aux9}). For this purpose we differentiate the two expressions for $\Phi^2$ in (\ref{Phi2}), (\ref{aux9}), simplify them by means of formulas (\ref{integrplushigh}), (\ref{integr11}), (\ref{integr3}), (\ref{aux5}), (\ref{curv3}), (\ref{aux18}) and compare the terms involving $\omega^1\wedge\omega^{\bar 1}$. This yields
$$
\Phi^2_0\equiv \frac{1}{2}\Theta^2_{21}\varphi^1+\frac{1}{2}\Theta^{\bar 2}_{{\bar 2}{\bar 1}}\varphi^{\bar 1}\quad(\mod\omega,\omega^1,\omega^{\bar 1},\theta^2,\theta^{\bar 2}).
$$
Thus, using (\ref{curv3}) we see that formula (\ref{aux9}) turns into
\begin{equation}
\begin{array}{l}
\displaystyle\Phi^2=\Theta^2_{21}\theta^2\wedge\omega^{\bar 1}+\Theta^{\bar 2}_{{\bar 2}{\bar 1}}\omega^1\wedge\theta^{\bar 2}+\frac{1}{2}\Theta^2_{21}\varphi^1\wedge\omega+\frac{1}{2}\Theta^{\bar 2}_{{\bar 2}{\bar 1}}\varphi^{\bar 1}\wedge\omega+\\
\vspace{-0.1cm}\\
\hspace{1cm}\Phi^2_{20}\theta^2\wedge\omega+\Phi^{\bar 2}_{{\bar 2}{\bar 0}}\theta^{\bar 2}\wedge\omega+\Phi^2_{10}\omega^1\wedge\omega+\Phi^{\bar 2}_{{\bar 1}{\bar 0}}\omega^{\bar 1}\wedge\omega,
\end{array}\label{curvPhi2}
\end{equation}
with
\begin{equation}
\Phi^2_0=\frac{1}{2}\Theta^2_{21}\varphi^1+\frac{1}{2}\Theta^{\bar 2}_{{\bar 2}{\bar 1}}\varphi^{\bar 1}+\Phi^2_{20}\theta^2+\Phi^{\bar 2}_{{\bar 2}{\bar 0}}\theta^{\bar 2}+\Phi^2_{10}\omega^1+\Phi^{\bar 2}_{{\bar 1}{\bar 0}}\omega^{\bar 1}.\label{aux15}
\end{equation} 

Next, (\ref{transform}), (\ref{Phi1}) imply
$$
\tilde\Phi^1\equiv \Phi^1+\frac{3r}{2}\omega^1\wedge\omega^{\bar 1}\quad(\mod\omega).
$$
By formula (\ref{aux15}) the expansion of $\Phi^1$ in (\ref{aux10}) does not involve $\omega^1\wedge\psi$, $\omega^{\bar 1}\wedge\psi$, hence
$$
\tilde\Phi^1_{1\bar 1}=\Phi^1_{1\bar 1}+\frac{3r}{2},
$$
which shows that $r$ is uniquely fixed by the requirement $\tilde\Phi^1_{1\bar 1}=0$. Thus, from now on we assume that the condition
\begin{equation}
\Phi^1_{1\bar 1}=0\label{curv4}
\end{equation}
is satisfied, hence in formula (\ref{transform}) one has $c=f=g=r=0$. Observe that $\varphi^1$ has now been canonically fixed.  Also, it is immediate from (\ref{aux10}), (\ref{aux15}) that $\Phi^1_{1\bar 1}=\Theta^2_{10}-\Phi^{\bar 2}_{{\bar 1}{\bar 0}}$, hence condition (\ref{curv4}) is equivalent to
$$
\Theta^2_{10}=\Phi^{\bar 2}_{{\bar 1}{\bar 0}}.
$$

We will now fix the function $s$. For this purpose we use the following form:
\begin{equation}
\Psi:=d\psi+\varphi^1\wedge\varphi^{\bar 1}+\varphi\wedge\psi\label{Psi}
\end{equation} 
(cf. equation (f) in (\ref{structure1})). Exterior differentiation of (\ref{integr2}) and substitution of (\ref{integrplushigh}), (\ref{integr11}) for $d\omega$, $d\omega^1$, respectively, yield
\begin{equation}
\begin{array}{l}
(\theta^{\bar 2}\wedge\varphi^1-\varphi^{\bar 1}\wedge\varphi^2+d\varphi^{\bar 1}+\omega^{\bar 1}\wedge\psi)\wedge\omega^1+\\
\vspace{-0.1cm}\\
(\theta^2\wedge\varphi^{\bar 1}-\varphi^1\wedge\varphi^{\bar 2}+d\varphi^1-\omega^1\wedge\psi)\wedge\omega^{\bar 1}+\\
\vspace{-0.1cm}\\
2(\varphi^1\wedge\varphi^{\bar 1}+\varphi\wedge\psi+d\psi)\wedge\omega=0.
\end{array}\label{exp4}
\end{equation}
Now, formulas (\ref{Phi1}), (\ref{aux10}), (\ref{exp4}) imply
\begin{equation}
d\psi\equiv-\varphi^1\wedge\varphi^{\bar 1}-\varphi\wedge\psi-\frac{1}{2}\Phi^{\bar 1}_{\bar 0}\wedge\omega^1+\frac{1}{2}\Phi^1_0\wedge\omega^{\bar 1} \quad (\mod \omega),\label{aux24}
\end{equation}
hence
\begin{equation}
\Psi=-\frac{1}{2}\Phi^{\bar 1}_{\bar 0}\wedge\omega^1+\frac{1}{2}\Phi^1_0\wedge\omega^{\bar 1}+\Psi_0\wedge\omega\label{aux11}
\end{equation}
for some real-valued 1-form $\Psi_0$. 

Before proceeding further, we analyze the expansion of $\Phi^1_0$ in (\ref{aux10}). For this purpose we differentiate the two expressions for $\Phi^1$ in (\ref{Phi1}), (\ref{aux10}), simplify them by means of formulas (\ref{integrplushigh}), (\ref{integr11}), (\ref{integr3}), (\ref{aux5}), (\ref{Phi2}), (\ref{aux9}), (\ref{curv3}), (\ref{aux18}), (\ref{aux15}), (\ref{curv4}), (\ref{aux24}) and compare the terms involving $\omega^1\wedge\omega^{\bar 1}$. This yields
$$
\Phi^1_0\equiv 0\quad(\mod\omega,\omega^1,\omega^{\bar 1},\theta^2,\theta^{\bar 2},\varphi^1,\varphi^{\bar 1},\psi).
$$
Thus, using (\ref{aux15}) (\ref{curv4}) we see that formula (\ref{aux10}) turns into
\begin{equation}
\begin{array}{l}
\displaystyle\Phi^1=\Theta^2_{20}\theta^2 \wedge\omega^{\bar 1}-\Phi^{\bar 2}_{{\bar 2}{\bar 0}}\omega^1\wedge\theta^{\bar 2}+\Phi^2_{20}\theta^2\wedge\omega^1-
\frac{1}{2}\Theta^2_{21}\omega^1\wedge\varphi^1-\\
\vspace{-0.1cm}\\
\displaystyle\hspace{1cm}\frac{1}{2}\Theta^{\bar 2}_{{\bar 2}{\bar 1}}\omega^1\wedge\varphi^{\bar 1}+{\bf P}_1\varphi^1\wedge\omega+{\bf P}_2\varphi^{\bar 1}\wedge\omega+{\bf P}_3\psi\wedge\omega+\\
\vspace{-0.1cm}\\
\hspace{1cm}\Phi^1_{20}\theta^2\wedge\omega+\Phi^1_{\bar 2 0}\theta^{\bar 2}\wedge\omega+\Phi^1_{10}\omega^1\wedge\omega+\Phi^1_{\bar 10}\omega^{\bar 1}\wedge\omega,
\end{array}\label{aux28}
\end{equation}with
\begin{equation}
\Phi^1_0={\bf P}_1\varphi^1+{\bf P}_2\varphi^{\bar 1}+{\bf P}_3\psi+\Phi^1_{20}\theta^2+\Phi^1_{\bar 2 0}\theta^{\bar 2}+\Phi^1_{10}\omega^1+\Phi^1_{\bar 10}\omega^{\bar 1}\label{Phi10}
\end{equation}
for some functions ${\bf P}_j$. In addition, our calculations show that all the functions ${\bf P}_j$ identically vanish whenever all the functions $\Theta^2_{21}$, $\Theta^2_{20}$, $\Phi^2_{20}$ do.

Next, (\ref{transform}), (\ref{Psi}) imply
$$
\tilde\Psi\equiv\Psi+s\omega^1\wedge\omega^{\bar 1}\quad(\mod\omega),
$$
hence
$$
\tilde\Psi_{1\bar 1}=\Psi_{1\bar 1}+s,
$$
which shows that $s$ is uniquely fixed by the requirement $\tilde\Psi_{1\bar 1}=0$ (observe that $\Psi_{1\bar 1}$ is real-valued). Thus, from now on we assume that the condition
\begin{equation}
\Psi_{1\bar 1}=0\label{curv5}
\end{equation}
is satisfied, therefore $\psi$ has now been canonically fixed. Also, it is immediate from (\ref{aux11}), (\ref{Phi10}) that $\Psi_{1\bar 1}=\Re\Phi^1_{10}$, hence condition (\ref{curv5}) is equivalent to
$$
\Re\Phi^1_{10}=0.
$$

We will now finalize our expansion for $\Psi$ by analyzing the form $\Psi_0$ in (\ref{aux11}). For this purpose we differentiate the two expressions for $\Psi$ in (\ref{Psi}), (\ref{aux11}), simplify them by means of formulas (\ref{integrplushigh}), (\ref{integr11}), (\ref{integr2}), (\ref{integr3}), (\ref{aux5}), (\ref{aux18}), (\ref{aux15}), (\ref{curv4}), (\ref{aux24}), (\ref{Phi10}), (\ref{curv5}) and compare the terms involving $\omega^1\wedge\omega^{\bar 1}$. This yields
$$
\Psi_0\equiv 0\quad(\mod\omega,\omega^1,\omega^{\bar 1},\theta^2,\theta^{\bar 2},\varphi^1,\varphi^{\bar 1},\psi).
$$
Thus, using (\ref{Phi10}), (\ref{curv5}) we see that formula (\ref{aux11}) turns into
\begin{equation}
\begin{array}{l}
\displaystyle\Psi=\frac{1}{2}\Phi^1_{20}\theta^2\wedge\omega^{\bar 1}+\frac{1}{2}\Phi^{\bar 1}_{\bar 2\bar 0}\omega^1\wedge\theta^{\bar 2}-\frac{1}{2}\Phi^{\bar 1}_{2\bar 0}\theta^2\wedge\omega^1+\frac{1}{2}\Phi^1_{\bar 2 0}\theta^{\bar 2}\wedge\omega^{\bar 1}+\\
\vspace{-0.1cm}\\
\displaystyle\hspace{0.9cm}\frac{1}{2}{\bf P}_{\bar 2}\omega^1\wedge\varphi^1+\frac{1}{2}{\bf P}_{\bar 1}\omega^1\wedge\varphi^{\bar 1}-\frac{1}{2}{\bf P}_{\bar 3}\omega^1\wedge\psi+\frac{1}{2}{\bf P}_1\varphi^1\wedge\omega^{\bar 1}+\\
\vspace{-0.1cm}\\
\displaystyle\hspace{0.9cm}\frac{1}{2}{\bf P}_2\varphi^{\bar 1}\wedge\omega^{\bar 1}+\frac{1}{2}{\bf P}_3\psi\wedge\omega^{\bar 1}+{\bf Q}_1\varphi^1\wedge\omega+{\bf Q}_{\bar 1}\varphi^{\bar 1}\wedge\omega+\\
\vspace{-0.1cm}\\
\displaystyle\hspace{0.9cm}{\bf Q}_3\psi\wedge\omega+\Psi_{20}\theta^2\wedge\omega+\Psi_{\bar 2 \bar 0}\theta^{\bar 2}\wedge\omega+\Psi_{10}\omega^1\wedge\omega+\Psi_{\bar 1 \bar 0}\omega^{\bar 1}\wedge\omega
\end{array}\label{aux48}
\end{equation}
with
$$
\Psi_0={\bf Q}_1\varphi^1+{\bf Q}_{\bar 1}\varphi^{\bar 1}+{\bf Q}_3\psi+\Psi_{20}\theta^2+\Psi_{\bar 2 \bar 0}\theta^{\bar 2}+\Psi_{10}\omega^1+\Psi_{\bar 1\bar 0}\omega^{\bar 1}
$$
for some functions ${\bf Q}_j$, where ${\bf Q}_3$ is $i\RR$-valued. In addition, our calculations show that all ${\bf Q}_j$ identically vanish whenever all the functions ${\bf P}_1$, ${\bf P}_2$, ${\bf P}_3$, $\Phi^1_{20}$ do.

Now that the forms $\theta^2$, $\varphi^1$, $\varphi^2$, $\psi$ have been uniquely determined by requirements (\ref{curv1}), (\ref{curv2}), (\ref{curv3}), (\ref{curv4}), (\ref{curv5}), they give rise to 1-forms defined on all of ${\mathcal P}^2$, and we denote these globally defined forms by the same respective symbols. Accordingly, the expansions for the 2-forms $\Theta^2$, $\Phi^1$, $\Phi^2$, $\Psi$ found in (\ref{curvv1}), (\ref{curvPhi2}), (\ref{aux28}), (\ref{aux48}) hold globally on ${\mathcal P}^2$.  

Further, let us denote by ${\mathcal P}_M$ the bundle ${\mathcal P}^2$ viewed as a principal $H$-bundle over $M$ and introduce a ${\mathfrak g}$-valued absolute parallelism $\omega_M$ on ${\mathcal P}_M$ by the formula
$$
\omega_M:=\left(
\begin{array}{ccccc}
\varphi^2 & \theta^2 & \omega^1&\omega& 0\\
\theta^{\bar 2} & \varphi^{\bar 2} &  \omega^{\bar 1}&0&-\omega\\
\varphi^{\bar 1} & \varphi^1 & 0&-\omega^{\bar 1}& -\omega^1\\
\psi & 0 & -\varphi^1&- \varphi^{\bar 2}& -\theta^2\\
0 & -\psi & -\varphi^{\bar 1}&-\theta^{\bar 2}& -\varphi^2\\
\end{array}\right)
$$
(cf. (\ref{mc})). Next, let ${\mathfrak h}\subset{\mathfrak g}$ be the Lie algebra of $H$. For every element $v\in{\mathfrak h}$ define $X_v$ to be the fundamental vector field on ${\mathcal P}_M$ arising from $v$, i.e.
$$
X_v({\bf\Theta}):=\frac{d}{dt}\Bigl(\exp(tv){\bf\Theta}\Bigr)\Bigl|_{t=0},\quad {\bf\Theta}\in{\mathcal P}_M.
$$
Using formulas (\ref{actionbyh}) it is straightforward to check that the parallelism $\omega_M$ has the following property:
\begin{equation}
\omega_M({\bf\Theta})(X_v({\bf\Theta})) = v\quad \hbox{for all $v\in{\mathfrak h}$ and ${\bf\Theta}\in{\mathcal P}_M$}.\label{ver}
\end{equation}

We are now ready to state the main theorem of the paper.

\begin{theorem}\label{main}\sl The CR-structures in the class ${\mathfrak C}_{2,1}$ are reducible to absolute parallelisms. Namely, for any $M,\tilde M\in{\mathfrak C}_{2,1}$ the following holds:

\noindent {\rm (i)} any CR-isomorphism $f:M\ra\tilde M$ can be uniquely lifted to a bundle isomorphism\linebreak $F:{\mathcal P}_{M}\ra{\mathcal P}_{\tilde M}$ satisfying
\begin{equation}
F^*\omega_{\tilde M}=\omega_{M},\label{equivf}
\end{equation}
{\rm (ii)} any diffeomorphism $F:{\mathcal P}_{M}\ra{\mathcal P}_{\tilde M}$ satisfying {\rm (\ref{equivf})} is a bundle isomorphism that is a lift of a CR-isomorphism $f:M\ra \tilde M$.        
\end{theorem}

\noindent{\bf Proof:} The existence of a lift $F$ in part (i) is immediate from the construction of the bundles and parallelisms. To see that such a lift is unique, let $F_1,F_2:{\mathcal P}_{M}\ra{\mathcal P}_{\tilde M}$ be two lifts of a CR-map $f:M\ra \tilde M$. Then ${\mathcal F}:=F_2^{-1}\circ F_1$ is a lift of the identity map of $M$ that preserves $\omega_M$. In particular, ${\mathcal F}$ preserves the 1-forms $\omega$, $\omega^1$ and $\varphi$. It then follows that ${\mathcal F}$ is the identity map, hence $F_1\equiv F_2$.

To obtain part (ii), suppose that $F:{\mathcal P}_{M}\ra{\mathcal P}_{\tilde M}$ is a diffeomorphism satisfying (\ref{equivf}). Then property (\ref{ver}) implies that $F$ is in fact a fiber bundle isomorphism i.e.~$F$ is a lift of a diffeomorphism $f: M\ra \tilde M$. We will now show that $f$ is a CR-map. Let $p\in M$. Since $F^*\tilde\omega=\omega$, for any covector $\tilde\theta$ at $f(p)$ annihilating $H_{f(p)}(\tilde M)$, the covector $f^*\tilde\theta$ at $p$ annihilates $H_{p}(M)$. Next, since $F^*\tilde\omega^1=\omega^1$, for any covector $\tilde\theta_0^1$ at $f(p)$ which is complex-linear on $H_{f(p)}(\tilde M)$ and satisfies $\tilde\theta_0^1(\tilde X)=0$ if $\tilde X-i\tilde J_{f(p)}\tilde X\in\ker{\mathcal L}_{\tilde M}(f(p))$, the covector $\theta_0^1:=f^*\tilde\theta_0^1$ is complex-linear on $H_{p}(M)$ and satisfies $\theta_0^1(X)=0$ if $X-iJ_{p}X\in\ker{\mathcal L}_{M}(p)$. Finally, since $F^*\tilde\theta^2=\theta^2$, for any covector $\tilde\theta_0^1$ as above and some covector $\tilde\theta_0^2$ at $f(p)$ such that $\tilde\theta_0^2$ is complex-linear on $H_{f(p)}(\tilde M)$ and the restrictions of $\tilde\theta_0^1$ and $\tilde\theta_0^2$ to $H_{f(p)}(\tilde M)$ form a basis of $H_{f(p)}^*(\tilde M)$, the covector $\theta_0^2:=f^*\tilde\theta_0^2$ is complex-linear on $H_p(M)$ and the restrictions of $\theta_0^1:=f^*\tilde\theta_0^1$, $\theta_0^2$ to $H_{p}(M)$ form a basis of $H_{p}^*(M)$. Thus, we have shown that $df(p)$ maps $H_{p}(M)$ into $H_{f(p)}(\tilde M)$ and is complex-linear on $H_p(M)$, i.e.~$f$ is a CR-map.\qed
\vspace{0.3cm}

Next, we define the curvature $\Omega_M$ of $\omega_M$ as the following ${\mathfrak g}$-valued 2-form:
$$
\Omega_M:=d\omega_M-\omega_M\wedge\omega_M.
$$
Then in terms of matrix elements identities (\ref{integrplushigh}), (\ref{integr11}), (\ref{integr2}) can be written as 
$$
(\Omega_M)^1_4=0,\quad (\Omega_M)^1_3=0,\quad \Re(\Omega_M)^1_1=0,
$$
respectively, and for the 2-forms $\Theta^2$, $\Phi^1$, $\Phi^2$, $\Psi$ defined in (\ref{Sigma}), (\ref{Phi2}), (\ref{Phi1}), (\ref{Psi}), we have
$$
\Theta^2=(\Omega_M)^1_2,\quad\Phi^1=(\Omega_M)^3_2,\quad \Phi^2=(\Omega_M)^1_1,\quad \Psi=(\Omega_M)^4_1.
$$
In the next section we give an application of the expansions of the curvature components obtained in (\ref{curvv1}), (\ref{curvPhi2}), (\ref{aux28}), (\ref{aux48}).

\section{Leading curvature terms and Cartan connections}\label{curvature}
\setcounter{equation}{0}

For applications it would be desirable to construct a Cartan connection on ${\mathcal P}_M$, not just an absolute parallelism. To recall the definition of Cartan connection, let $R$ be a Lie group with Lie algebra ${\mathfrak r}$ and $S$ a closed subgroup of $R$ with Lie algebra ${\mathfrak s}\subset{\mathfrak r}$ acting by diffeomorphisms on a manifold ${\mathcal P}$ such that $\dim{\mathcal P}= \dim R$. As before, for every element $v\in{\mathfrak s}$ denote by $X_v$ the fundamental vector field on ${\mathcal P}$ arising from $v$. A Cartan connection of type $S\hspace{-0.1cm}\setminus\hspace{-0.1cm}R$ on the manifold ${\mathcal P}$ is an ${\mathfrak r}$-valued absolute parallelism $\pi$ satisfying
\begin{equation}
\begin{array}{ll}
{\rm (i)} & \hbox{$\pi(x)(X_v(x))=v$ for all $v\in{\mathfrak s}$ and $x\in{\mathcal P}$, and}\\
\vspace{-0.1cm}\\
{\rm (ii)} &\hbox{$L^*_s\pi=\hbox{Ad}_{S,{\mathfrak r}}(s)\pi$ for all $s\in S$,}
\end{array}\label{cartanconnec}
\end{equation}
where $L_s$ denotes the (left) action by an element $s$ on ${\mathcal P}$ and $\hbox{Ad}_{S,{\mathfrak r}}$ is the adjoint representation of $S$. For example, the Maurer-Cartan form $\omega_R^{\hbox{\tiny MC}}$ on $R$ is a Cartan connection of type $S\hspace{-0.1cm}\setminus\hspace{-0.1cm}R$ for any subgroup $S$. 

As we will see below, the expansions of components of the curvature form obtained in (\ref{curvPhi2}), (\ref{aux28}), (\ref{aux48}) can be used for determining obstructions for the parallelism $\omega_M$ to be a Cartan connection. Specifically, for the curvature components $\Phi^1$, $\Phi^2$, $\Psi$ we consider the leading terms, i.e.~the terms that involve either $\theta^2\wedge\omega^{\bar 1}$ or $\omega^1\wedge\theta^{\bar 2}$. From (\ref{curvPhi2}), (\ref{aux28}), (\ref{aux48}) we can write these terms explicitly as follows:
\begin{equation}
\begin{array}{l}
\Phi^1=\Theta^2_{20}\theta^2\wedge\omega^{\bar 1}-\Phi^{\bar 2}_{\bar 2\bar 0}\omega^1\wedge\theta^{\bar 2}+\dots,\\
\vspace{-0.1cm}\\
\Phi^2=\Theta^2_{21}\theta^2\wedge\omega^{\bar 1}+\Theta^{\bar 2}_{{\bar 2}{\bar 1}}\omega^1\wedge\theta^{\bar 2}+\dots,\\
\vspace{-0.1cm}\\
\displaystyle\Psi=\frac{1}{2}\Phi^1_{20}\theta^2\wedge\omega^{\bar 1}+\frac{1}{2}\Phi^{\bar 1}_{\bar 2\bar 0}\omega^1\wedge\theta^{\bar 2}+\dots.
\end{array}\label{aux88888}
\end{equation}
The result of this section is the following theorem.

\begin{theorem}\label{carconnec}\sl The absolute parallelism $\omega_M$ is a Cartan connection of type $H\hspace{-0.1cm}\setminus\hspace{-0.1cm}G$ if and only if all leading curvature terms identically vanish.
\end{theorem}

\noindent{\bf Proof:} We will first obtain the sufficiency implication. As was observed in (\ref{ver}), part (i) of (\ref{cartanconnec}) holds for $\omega_M$ with $R=G$, $S=H$, thus we only need to verify the identity
\begin{equation}
L^*_h\omega_M=\hbox{Ad}_{H,{\mathfrak g}}(h)\omega_M\label{cartanid}
\end{equation}
for all $h\in H$. Recall that $H=H^1\ltimes H^2$ (see (\ref{subgroups})). Fix $h\in H$ and write it as $h=h_1h_2$ with $h_j\in H^j$. 

First, let 
$$
\hat\omega_M:=\hbox{Ad}_{H,{\mathfrak g}}(h_2)\omega_M= 
\left(
\begin{array}{ccccc}
\hat\varphi^2 & \hat\theta^2 & \hat\omega^1&\hat\omega& 0\\
\hat\theta^{\bar 2} & \hat\varphi^{\bar 2} &  \hat\omega^{\bar 1}&0&-\hat\omega\\
\hat\varphi^{\bar 1} & \hat\varphi^1 & 0&-\hat\omega^{\bar 1}& -\hat\omega^1\\
\hat\psi & 0 & -\hat\varphi^1&- \hat\varphi^{\bar 2}& -\hat\theta^2\\
0 & -\hat\psi & -\hat\varphi^{\bar 1}&-\hat\theta^{\bar 2}& -\hat\varphi^2\\
\end{array}\right)
$$
for some 1-forms $\hat\omega$, $\hat\omega^1$, $\hat\theta^2$, $\hat\varphi^1$, $\hat\varphi^2$, $\hat\psi$. Writing
$$
h_2=\left(
\begin{array}{ccccc}
1 & 0 & 0&0& 0\\
0 & 1 & 0&0& 0\\
B & \bar B & 1&0& 0\\
\Lambda-|B|^2/2 & -\bar B^2/2 & -\bar B&1& 0\\
-B^2/2 & -\Lambda-|B|^2/2 & - B&0&1\\
\end{array}
\right)
$$
for some $B\in\CC$ and $\Lambda\in i\RR$ (see (\ref{subgroups})), we calculate
\begin{equation}
\begin{array}{l}
\displaystyle\hat\omega=\omega,\quad\hat\omega^1=\omega^1+{\bar B}\omega,\quad\hat\theta^2=\theta^2-{\bar B}\omega^1-\frac{{\bar B}^2}{2}\omega,\\
\vspace{-0.1cm}\\
\displaystyle\hat\varphi^1=\varphi^1-\left(\Lambda+\frac{|B|^2}{2}\right)\omega^1-\frac{{\bar B}^2}{2}\omega^{\bar 1}+B\theta^2-\Lambda{\bar B}\omega+{\bar B}\varphi^{\bar 2},\\
\vspace{-0.1cm}\\
\displaystyle\hat\varphi^2=\varphi^2-B\omega^1-\left(\Lambda+\frac{|B|^2}{2}\right)\omega,\\
\vspace{-0.1cm}\\
\displaystyle\hat\psi=\psi-\Lambda B\omega^1-\Lambda{\bar B}\omega^{\bar 1}+\frac{B^2}{2}\theta^2-\frac{{\bar B}^2}{2}\theta^{\bar 2}-\Lambda^2\omega+\\
\vspace{-0.3cm}\\
\displaystyle\hspace{0.8cm}B\varphi^1-{\bar B}\varphi^{\bar 1}+\left(\Lambda-\frac{|B|^2}{2}\right)\varphi^2+\left(\Lambda+\frac{|B|^2}{2}\right)\varphi^{\bar 2}.
\end{array}\label{transformhats}
\end{equation}
Analogously to formulas (\ref{Sigma}), (\ref{Phi2}), (\ref{Phi1}), (\ref{Psi}) we now define
\begin{equation}
\begin{array}{l}
\hat\Theta^2:=d\hat\theta^2+\hat\theta^2\wedge(\hat\varphi^2-\hat\varphi^{\bar 2})-\hat\omega^1\wedge\hat\varphi^1,\\
\vspace{-0.1cm}\\
\displaystyle\hat\Phi^1:=d\hat\varphi^1+\hat\theta^2\wedge\hat\varphi^{\bar 1}-\hat\omega^1\wedge\hat\psi-\hat\varphi^1\wedge\hat\varphi^{\bar 2},\\
\vspace{-0.1cm}\\
\hat\Phi^2:=d\hat\varphi^2-\hat\theta^2\wedge\hat\theta^{\bar 2}-\hat\omega^1\wedge\hat\varphi^{\bar 1}-\hat\omega\wedge\hat\psi,\\
\vspace{-0.1cm}\\
\hat\Psi:=d\hat\psi+\hat\varphi^1\wedge\hat\varphi^{\bar 1}+(\hat\varphi^2+\hat\varphi^{\bar 2})\wedge\hat\psi.
\end{array}\label{deflargehats}
\end{equation}
Observe that the 2-forms introduced in (\ref{deflargehats}) are components of $\hat\Omega_M:=\hbox{Ad}_{H,{\mathfrak g}}(h_2)\Omega_M$. Using identities (\ref{integrplushigh}), (\ref{integr11}), (\ref{transformhats}) we then see
\begin{equation}
\begin{array}{l}
\displaystyle\hat\Theta^2=\Theta^2,\quad \displaystyle\hat\Phi^1=\Phi^1+B\Theta^2-{\bar B}\Phi^2,\quad \hat\Phi^2=\Phi^2,\\
\vspace{-0.1cm}\\
\displaystyle\hat\Psi=\Psi+\frac{B^2}{2}\Theta^2-\frac{{\bar B}^2}{2}\Theta^{\bar 2}+B\Phi^1-{\bar B}\Phi^{\bar 1}-|B|^2\Phi^2.
\end{array}\label{aux888}
\end{equation}

Further, since $\Theta^2_{21}=\Theta^2_{20}=\Phi^1_{20}=\Phi^2_{20}=0$ (see (\ref{aux88888})), expansions (\ref{curvv1}), (\ref{curvPhi2}), (\ref{aux28}), (\ref{aux48}) become
\begin{equation}
\begin{array}{l}
\displaystyle \Theta^2 =\Theta^2_{10}\omega^1\wedge\omega
+\Theta^2_{\bar 10} \omega^{\bar 1}\wedge\omega,\\
\vspace{-0.1cm}\\
\displaystyle\Phi^1=\Phi^1_{\bar 2 0}\theta^{\bar 2}\wedge\omega+\Phi^1_{10}\omega^1\wedge\omega+\Phi^1_{\bar 10}\omega^{\bar 1}\wedge\omega,\\
\vspace{-0.1cm}\\
\displaystyle\Phi^2=\Phi^2_{10}\omega^1\wedge\omega+\Phi^{\bar 2}_{{\bar 1}{\bar 0}}\omega^{\bar 1}\wedge\omega,\\
\vspace{-0.1cm}\\
\displaystyle\Psi=-\frac{1}{2}\Phi^{\bar 1}_{2\bar 0}\theta^2\wedge\omega^1+\frac{1}{2}\Phi^1_{\bar 2 0}\theta^{\bar 2}\wedge\omega^{\bar 1}+\Psi_{20}\theta^2\wedge\omega+\\
\vspace{-0.3cm}\\
\hspace{1.2cm}\Psi_{\bar 2 \bar 0}\theta^{\bar 2}\wedge\omega+\Psi_{10}\omega^1\wedge\omega+\Psi_{\bar 1 \bar 0}\omega^{\bar 1}\wedge\omega.
\end{array}\label{aux444}
\end{equation}
Now formulas (\ref{aux888}), (\ref{aux444}) yield
\begin{equation}
\makebox[250pt]{$\begin{array}{l}
\displaystyle \hat\Theta^2 =\Theta^2_{10}\omega^1\wedge\omega
+\Theta^2_{\bar 10} \omega^{\bar 1}\wedge\omega,\\
\vspace{-0.1cm}\\
\displaystyle\hat\Phi^1=\Phi^1_{\bar 2 0}\theta^{\bar 1}\wedge\omega+(\Phi^1_{10}+B\Theta^2_{10}-{\bar B}\Phi^2_{10})\omega^1\wedge\omega+(\Phi^1_{\bar 1 0}+B\Theta^2_{\bar 1 0}-{\bar B}\Phi^{\bar 2}_{\bar 1\bar 0})\omega^{\bar 1}\wedge\omega,\\
\vspace{-0.1cm}\\
\displaystyle\hat\Phi^2=\Phi^2_{10}\omega^1\wedge\omega+\Phi^{\bar 2}_{{\bar 1}{\bar 0}}\omega^{\bar 1}\wedge\omega,\\
\vspace{-0.1cm}\\
\displaystyle\hat\Psi=-\frac{1}{2}\Phi^{\bar 1}_{2\bar 0}\theta^2\wedge\omega^1+\frac{1}{2}\Phi^1_{\bar 2 0}\theta^{\bar 2}\wedge\omega^{\bar 1}+(\Psi_{20}+{\bar B}\Phi^{\bar 1}_{2\bar 0})\theta^2\wedge\omega+(\Psi_{\bar 2\bar 0}+B\Phi^1_{\bar 2 0})\theta^{\bar 2}\wedge\omega+\\
\vspace{-0.3cm}\\
\displaystyle\hspace{1cm}\left(\Psi_{10}+\frac{B^2}{2}\Theta^2_{10}+\frac{{\bar B}^2}{2}\Theta^{\bar 2}_{1\bar 0}+B\Phi^1_{10}+{\bar B}\Phi^{\bar 1}_{1\bar 0}-|B|^2\Phi^2_{10}\right)\omega^1\wedge\omega+\\
\vspace{-0.3cm}\\
\displaystyle\hspace{1cm}\left(\Psi_{\bar 1\bar 0}+\frac{{\bar B}^2}{2}\Theta^{\bar 2}_{\bar 1 \bar 0}+\frac{B^2}{2}\Theta^2_{\bar 1 0}+{\bar B}\Phi^{\bar 1}_{\bar 1\bar 0}+B\Phi^1_{\bar 1 0}-|B|^2\Phi^{\bar 2}_{\bar 1\bar 0}\right)\omega^{\bar 1}\wedge\omega.
\end{array}$}\label{aux8888}
\end{equation}

Next, let 
$$
\hbox{Ad}_{H,{\mathfrak g}}(h)\omega_M=\hbox{Ad}_{H,{\mathfrak g}}(h_1)\hat\omega_M=
\left(
\begin{array}{ccccc}
\check \varphi^2 &\check\theta^2 & \check\omega^1&\check\omega& 0\\
\check\theta^{\bar 2} & \check\varphi^{\bar 2} & \check\omega^{\bar 1}&0&-\check\omega\\
\check\varphi^{\bar 1} & \check\varphi^1 & 0&-\check\omega^{\bar 1}& -\check\omega^1\\
\check\psi & 0 & -\check\varphi^1&- \check\varphi^{\bar 2}& -\check\theta^2\\
0 & -\check\psi & -\check\varphi^{\bar 1}&-\check\theta^{\bar 2}& -\check\varphi^2\\
\end{array}\right)
$$
for some 1-forms $\check\omega$, $\check\omega^1$, $\check\theta^2$, $\check\varphi^1$, $\check\varphi^2$, $\check\psi$. Writing
$$
h_1=\left(
\begin{array}{ccccc}
A & 0 & 0&0& 0\\
0 & \bar A & 0&0& 0\\
0 & 0 & 1&0& 0\\
0 & 0 & 0&\bar A^{-1}& 0\\
0 & 0 & 0&0&A^{-1}\\
\end{array}
\right)
$$
for some $A\in\CC^*$ (see (\ref{subgroups})), we calculate
\begin{equation}
\begin{array}{lll}
\displaystyle\check\omega=|A|^2\hat\omega,& \displaystyle\check\omega^1=A\hat\omega^1,&\displaystyle\check\theta^2=\frac{A}{\bar A}\hat\theta^2,\\
\vspace{-0.1cm}\\
\displaystyle\check\varphi^1=\frac{1}{\bar A}\hat\varphi^1,&\displaystyle\check\varphi^2=\hat\varphi^2,&\displaystyle\check\psi=\frac{1}{|A|^2}\hat\psi.
\end{array}\label{h1hats}
\end{equation}
From formulas (\ref{actionbyh}), (\ref{transformhats}), (\ref{h1hats}) we now observe
$$
\check\omega=L^*_h\omega,\quad \check\omega^1=L^*_h\omega^1,\quad \check\varphi=L^*_h\varphi,
$$
where $\check\varphi:=2\Re\check\varphi^2$. Thus, in order to obtain (\ref{cartanid}), one needs to show that the 1-forms with checks satisfy identities (\ref{integrplushigh}), (\ref{integr11}), (\ref{integr2}) as well as curvature conditions (\ref{curv1}), (\ref{curv2}), (\ref{curv3}), (\ref{curv4}), (\ref{curv5}). Straightforward calculations yield that this is indeed the case for identities (\ref{integrplushigh}), (\ref{integr11}), (\ref{integr2}) irrespectively of the vanishing of the leading curvature terms.

To deal with conditions (\ref{curv1}), (\ref{curv2}), (\ref{curv3}), (\ref{curv4}), (\ref{curv5}), we define, analogously to formulas (\ref{Sigma}), (\ref{Phi2}), (\ref{Phi1}), (\ref{Psi})
\begin{equation}
\begin{array}{l}
\check\Theta^2:=d\check\theta^2+\check\theta^2\wedge(\check\varphi^2-\check\varphi^{\bar 2})-\check\omega^1\wedge\check\varphi^1,\\
\vspace{-0.1cm}\\
\displaystyle\check\Phi^1:=d\check\varphi^1+\check\theta^2\wedge\check\varphi^{\bar 1}-\check\omega^1\wedge\check\psi-\check\varphi^1\wedge\check\varphi^{\bar 2},\\
\vspace{-0.1cm}\\
\check\Phi^2:=d\check\varphi^2-\check\theta^2\wedge\check\theta^{\bar 2}-\check\omega^1\wedge\check\varphi^{\bar 1}-\check\omega\wedge\check\psi,\\
\vspace{-0.1cm}\\
\check\Psi:=d\check\psi+\check\varphi^1\wedge\check\varphi^{\bar 1}+\check\varphi\wedge\check\psi.
\end{array}\label{deflargechecks}
\end{equation}
Observe that the 2-forms introduced in (\ref{deflargechecks}) are components of $\hbox{Ad}_{H,{\mathfrak g}}(h)\Omega_M=\hbox{Ad}_{H,{\mathfrak g}}(h_1)\hat\Omega_M$. Our aim is to verify the identities
\begin{equation}
\check\Theta^2_{2\bar 1}=\check\Theta^2_{1\bar 1}=\check\Phi^1_{1\bar 1}=\check\Phi^2_{1\bar 1}=\check\Psi_{1\bar 1}=0.\label{idsss}
\end{equation}
By (\ref{h1hats}) we have 
\begin{equation}
\displaystyle\check\Theta^2=\frac{A}{\bar A}\hat\Theta^2,\quad
\displaystyle\check\Phi^1=\frac{1}{\bar A}\hat\Phi^1,\quad\displaystyle\check\Phi^2=\hat\Phi^2,\quad
\displaystyle\check\Psi=\frac{1}{|A|^2}\hat\Psi.\label{h1curvtransform}
\end{equation}
Identities (\ref{idsss}) are now an immediate consequence of (\ref{transformhats}), (\ref{aux8888}), (\ref{h1hats}), (\ref{h1curvtransform}).

Next, we obtain the necessity implication. Since $\omega_M$ is a Cartan connection, for every element $h\in H^2$ we have, in particular, $\hat\Phi_{1\bar1}=0$ and $\hat\Psi_{1\bar 1}=0$, where $\hat\Phi^1$, $\hat\Psi$ are defined in (\ref{deflargehats}). Using formulas (\ref{curvv1}), (\ref{curvPhi2}), (\ref{aux28}), (\ref{transformhats}), (\ref{aux888}) we calculate
$$
\hat\Phi^1_{1\bar 1}={\bar B}\Theta^2_{20}-B\Phi^{\bar 2}_{\bar 2\bar 0}-\frac{3{\bar B}^2}{4}\Theta^2_{21}-\frac{1}{2}\left(\Lambda+\frac{3|B|^2}{2}\right)\Theta^{\bar 2}_{\bar 2\bar 1}.
$$
The vanishing of $\hat\Phi^1_{1\bar 1}$ for any $B\in\CC$ and $\Lambda\in i\RR$ then implies $\Theta^2_{21}=\Theta^2_{20}=\Phi^2_{20}=0$. Therefore, expansion (\ref{aux48}) becomes
\begin{equation}
\begin{array}{l}
\displaystyle\Psi=\frac{1}{2}\Phi^1_{20}\theta^2\wedge\omega^{\bar 1}+\frac{1}{2}\Phi^{\bar 1}_{\bar 2\bar 0}\omega^1\wedge\theta^{\bar 2}-\frac{1}{2}\Phi^{\bar 1}_{2\bar 0}\theta^2\wedge\omega^1+\frac{1}{2}\Phi^1_{\bar 2 0}\theta^{\bar 2}\wedge\omega^{\bar 1}+\\
\vspace{-0.1cm}\\
\displaystyle\hspace{0.9cm}{\bf Q}_1\varphi^1\wedge\omega+{\bf Q}_{\bar 1}\varphi^{\bar 1}\wedge\omega+{\bf Q}_3\psi\wedge\omega+\Psi_{20}\theta^2\wedge\omega+\Psi_{\bar 2 \bar 0}\theta^{\bar 2}\wedge\omega+\\
\vspace{-0.1cm}\\
\displaystyle\hspace{0.9cm}\Psi_{10}\omega^1\wedge\omega+\Psi_{\bar 1 \bar 0}\omega^{\bar 1}\wedge\omega,
\end{array}\label{newPsiexp}
\end{equation}
and using formulas (\ref{curvv1}), (\ref{curvPhi2}), (\ref{aux28}), (\ref{transformhats}), (\ref{aux888}), (\ref{newPsiexp}) we calculate
$$
\hat\Psi_{1\bar 1}=\frac{{\bar B}}{2}\Phi^1_{20}+\frac{B}{2}\Phi^{\bar 1}_{\bar 2\bar 0}.
$$
The vanishing of $\hat\Psi_{1\bar 1}$ for any $B\in\CC$ then implies $\Phi^1_{20}=0$, which completes the proof of the theorem.\qed
\vspace{0.3cm}

We will now demonstrate that leading curvature terms can indeed occur. Namely, we give examples of hypersurfaces in $\CC^3$ for which the coefficient $\Theta^2_{21}$ does not identically vanish (cf. (\ref{aux88888})). By Theorem \ref{carconnec}, for any such hypersurface the absolute parallelism $\omega_M$ is not a Cartan connection.

Let $M$ be the hypersurface in $\CC^3$ given by the equation
$$
z_3+{\bar z}_3=\rho(z_1+{\bar z}_1,z_2+{\bar z}_2),
$$
where $\rho(t_1,t_2)$ is a smooth real-valued function on a domain in $\RR^2$. Assume that $\rho_{11}$ is everywhere positive and $\rho$ satisfies the Monge-Amp\`ere equation
\begin{equation}
\rho_{11}\rho_{22}-\rho_{12}^2\equiv 0,\label{mongeampere}
\end{equation}
where indices indicate partial derivatives of $\rho$, with index 1 corresponding to $t_1$ and index 2 to $t_2$. Then $M$ is uniformly Levi degenerate of rank 1.

From now on we assume that all functions of the variables $t_1$, $t_2$ are calculated for $t_1=z_1+{\bar z}_1$, $t_2=z_2+{\bar z}_2$. Using this convention and setting
\begin{equation}
\mu:=\rho_1dz_1+\rho_2dz_2-dz_3,\quad \eta^1:=\rho_{11}dz_1+\rho_{12}dz_2,\quad \eta^2:=dz_2,\label{aux45}
\end{equation}
we see
\begin{equation}
\begin{array}{l}
\displaystyle d\mu=-\frac{1}{\rho_{11}}\eta^1\wedge\eta^{\bar 1},\\
\vspace{-0.1cm}\\
\displaystyle d\eta^1=-\left(\frac{\rho_{12}}{\rho_{11}}\right)_{\hspace{-0.1cm}1}\eta^2\wedge\eta^{\bar 1}+\left[\frac{\rho_{111}}{\rho_{11}^2}\eta^{\bar 1}+\left(\frac{\rho_{12}}{\rho_{11}}\right)_{\hspace{-0.1cm}1}\eta^{\bar 2}\right]\wedge\eta^1.
\end{array}\label{aux7777}
\end{equation}
To obtain the second equation in (\ref{aux7777}) we use the identity
$$
\rho_{111}\left(\frac{\rho_{12}}{\rho_{11}}\right)^2-2\rho_{112}\left(\frac{\rho_{12}}{\rho_{11}}\right)+\rho_{122}=0,
$$
which is a consequence of (\ref{mongeampere}). Thus, $M$ is 2-nondegenerate if and only if the function $S:=(\rho_{12}/\rho_{11})_{1}$ vanishes nowhere, and we assume from now on that this condition is satisfied.

For the form $\omega$ on ${\mathcal P}^1$ we then have $\omega(u\mu)=u\mu^*$, therefore from the first equation in (\ref{aux7777}) one observes
\begin{equation}
d\omega= -\frac{u}{\rho_{1\bar 1}^*}\eta^{1*}\wedge\eta^{\bar 1*}-\omega\wedge\frac{du}{u},\label{domeganew}
\end{equation}
where $u>0$ is the fiber coordinate and asterisks indicate pull-backs from $M$ to ${\mathcal P}^1$. Identity (\ref{domeganew}) shows that setting
\begin{equation}
\nu:=\sqrt{\frac{u}{\rho_{11}^*}}\eta^{1*},\label{aux43}
\end{equation}
one can parametrize the fibers of ${\mathcal P}^2$ as
$$
\begin{array}{l}
\displaystyle\theta^1=a\nu+{\bar b}\omega,\\
\vspace{-0.1cm}\\
\displaystyle\phi=\frac{du}{u}-ab\nu-{\bar a}{\bar b}{\bar \nu}+\lambda\omega,
\end{array}\label{aux459}
$$
with $|a|=1$, $b\in\CC$, $\lambda\in i\RR$ (see (\ref{g1structure})). We then have
\begin{equation}
\begin{array}{l}
\displaystyle\omega^1=a\nu^*+{\bar b}\omega,\\
\vspace{-0.1cm}\\
\displaystyle\varphi=\left(\frac{du}{u}\right)^*-b\omega^1-{\bar b}\omega^{\bar 1}+\lambda\omega,
\end{array}\label{aux41}
\end{equation}
where asterisks indicate pull-backs from ${\mathcal P}^1$ to ${\mathcal P}^2$ and the pull-back of $\omega$ is denoted by the same symbol (cf. the notation of Section \ref{construction}).

Set
\begin{equation}
\theta^2:=-a^2S^{**}\eta^{2**}\label{aux84}
\end{equation}
where double asterisks indicate pull-backs from $M$ to ${\mathcal P}^2$. Identities (\ref{integrplushigh}), (\ref{aux45}), (\ref{aux7777}),  (\ref{aux43}), (\ref{aux41}), (\ref{aux84}) imply  
\begin{equation}
\begin{array}{l}
\displaystyle d\omega^1=\theta^2\wedge\omega^{\bar 1}-\omega^1\wedge\left[\frac{da}{a}+\left(\frac{du}{2u}\right)^*+\frac{\bar a^2}{2}\theta^2-\frac{a^2}{2}\theta^{\bar 2}+\right.\\
\vspace{-0.1cm}\\
\displaystyle\hspace{7cm}\left.\left({\bar b}+\frac{a\rho_{111}}{2\sqrt{u^*\rho_{11}^{3**}}}\right)\omega^{\bar 1}\right]-\omega\wedge(d{\bar b}+\sigma),
\end{array}\label{aux454}
\end{equation}
where $\sigma$ is a 1-form vanishing for $b=0$. Set
\begin{equation}
\begin{array}{l}
\displaystyle\varphi^1:=d{\bar b}+\frac{\lambda}{2}\omega^1+\sigma,\\
\vspace{-0.1cm}\\
\displaystyle\varphi^2:=\frac{da}{a}+\left(\frac{du}{2u}\right)^*+\frac{\bar a^2}{2}\theta^2-\frac{a^2}{2}\theta^{\bar 2}-\\
\vspace{-0.1cm}\\
\displaystyle\hspace{5cm}\left(2b+\frac{{\bar a}\rho_{111}}{2\sqrt{u^*\rho_{11}^{3**}}}\right)\omega^1+\left(\bar b+\frac{a\rho_{111}}{2\sqrt{u^*\rho_{11}^{3**}}}\right)\omega^{\bar 1}+\frac{\lambda}{2}\omega.
\end{array}\label{787}
\end{equation}
It follows from (\ref{aux41}), (\ref{aux454}), (\ref{787}) that the forms $\theta^2$, $\varphi^1$, $\varphi^2$ satisfy
\begin{equation}
\begin{array}{l}
\displaystyle\Re\varphi^2=\frac{\varphi}{2},\\
\vspace{-0.1cm}\\
\displaystyle d\omega^1=\theta^2\wedge\omega^{\bar 1}-\omega^1\wedge\varphi^2-\omega\wedge\varphi^1
\end{array}\label{aux111}
\end{equation}  
(cf. (\ref{cond1}), (\ref{integr11})).

For future arguments we will need expressions of $(du)^*$, $da$, $db$ in terms of $\omega$, $\omega^1$, $\omega^{\bar 1}$,  $\theta^2$, $\theta^{\bar 2}$, $\varphi^1$, $\varphi^{\bar 1}$, $\varphi^2$, $\varphi^{\bar 2}$. Identities (\ref{aux41}), (\ref{787}), (\ref{aux111}) yield
\begin{equation}
\begin{array}{l}
\displaystyle(du)^*=u^*(b\omega^1+{\bar b}\omega^{\bar 1}-\lambda\omega+\varphi^2+\varphi^{\bar 2}),\\
\vspace{-0.1cm}\\
\displaystyle da=a\left[-\frac{\bar a^2}{2}\theta^2+\frac{a^2}{2}\theta^{\bar 2}+\left(\frac{3b}{2}+\frac{{\bar a}\rho_{111}}{2\sqrt{u^*\rho_{11}^{3**}}}\right)\omega^1-\right.\\
\vspace{-0.3cm}\\
\displaystyle\hspace{7cm}\left.\left(\frac{3\bar b}{2}+\frac{a\rho_{111}}{2\sqrt{u^*\rho_{11}^{3**}}}\right)\omega^{\bar 1}+\frac{1}{2}\varphi^2-\frac{1}{2}\varphi^{\bar 2}\right],\\
\vspace{-0.1cm}\\
\displaystyle db=\frac{\lambda}{2}\omega^{\bar 1}-{\bar\sigma}+\varphi^{\bar 1}.
\end{array}\label{differentials}
\end{equation}  
Observe that for $u^*=1$, $a=1$, $b=0$, $\lambda=0$ formulas (\ref{differentials}) simplify as
\begin{equation}
\begin{array}{l}
\displaystyle(du)^*=\varphi^2+\varphi^{\bar 2},\\
\vspace{-0.1cm}\\
\displaystyle da=-\frac{1}{2}\theta^2+\frac{1}{2}\theta^{\bar 2}+\frac{\rho_{111}}{2\sqrt{\rho_{11}^{3**}}}\omega^1-\frac{\rho_{111}}{2\sqrt{\rho_{11}^{3**}}}\omega^{\bar 1}+\frac{1}{2}\varphi^2-\frac{1}{2}\varphi^{\bar 2},\\
\vspace{-0.1cm}\\
\displaystyle db=\varphi^{\bar 1}.
\end{array}\label{differentialsat0}
\end{equation}

Next, according to the normalization procedure described in Section \ref{construction}, we introduce\linebreak $\tilde\theta^2:=\theta^2-c\omega^1$ and $\tilde\varphi^2:=\varphi^2+{\bar c}\omega^1-c\omega^{\bar 1}$, where $c$ is chosen so that the expansion of $d\tilde\theta^2$ does not involve $\tilde\theta^2\wedge\omega^{\bar 1}$. By formula (\ref{choiceofc}), the function $c$ is given by
\begin{equation}
c=\frac{1}{3}\Theta^2_{2\bar 1},\label{paramc}
\end{equation}
where $\Theta^2_{2\bar 1}$ is the coefficient at the wedge product $\theta^2\wedge\omega^{\bar 1}$  in the expansion of $d\theta^2$ (see (\ref{Sigma})). From (\ref{aux45}), (\ref{aux43}), (\ref{aux41}), (\ref{aux84}), (\ref{differentials}) we find the following formula for this coefficient:
\begin{equation} 
\Theta^2_{2\bar 1}=-\frac{aS_{1}^{**}}{\sqrt{u\rho_{11}^{**}}S^{**}}+{3\bar b}+\frac{a\rho_{111}}{\sqrt{u^*\rho_{11}^{3**}}}.\label{aux87}
\end{equation} 
Further, by (\ref{integr11}), (\ref{Sigma}) one has
\begin{equation}
\displaystyle\tilde\Theta^2\equiv\Theta^2-dc\wedge\omega^1+2{\bar c}\theta^2\wedge\omega^1-3c\theta^2\wedge\omega^{\bar 1}\,\,(\mod\hbox{terms not involving $\theta^2$}).\label{tildetheta2spec}
\end{equation}

From now on we restrict our calculations to the section $\gamma_0$ of ${\mathcal P}^2$ given by $u^*=1$, $a=1$, $b=0$, $\lambda=0$. Using formulas (\ref{theta2}), (\ref{normal8}), (\ref{aux45}), (\ref{aux43}), (\ref{aux41}), (\ref{aux84}), (\ref{differentialsat0}), (\ref{paramc}), (\ref{aux87}), (\ref{tildetheta2spec}) we then obtain
\begin{equation}
\begin{array}{l}
\displaystyle\tilde\Theta^2_{21}=\Theta^2_{21}+\frac{1}{3S^{**}}\left[\frac{\rho_{12}^{**}}{\rho_{11}^{**}}\left(\frac{S_{1}^{**}}{\sqrt{\rho_{11}^{**}}S^{**}}\right)_{\hspace{-0.1cm}1}-\left(\frac{S_{1}^{**}}{\sqrt{\rho_{11}^{**}}S^{**}}\right)_{\hspace{-0.1cm}2}\right]-\\
\vspace{-0.1cm}\\
\displaystyle\hspace{4cm} \frac{1}{3S^{**}}\left[\frac{\rho_{12}^{**}}{\rho_{1 1}^{**}}\left(\frac{\rho_{111}}{\sqrt{\rho_{11}^{3**}}}\right)_{\hspace{-0.1cm}1}-\left(\frac{\rho_{111}}{\sqrt{\rho_{11}^{3**}}}\right)_{\hspace{-0.1cm}2}\right]-\frac{5S_{1}^{**}}{6\sqrt{\rho_{11}^{**}}\,\,{S^{**}}}+\frac{5\rho_{111}}{6\sqrt{\rho_{11}^{3**}}}.
\end{array}\label{finalxv}
\end{equation}
Notice that $\tilde\Theta^2_{21}$ in formula (\ref{finalxv}) is in fact the final value of this curvature coefficient since transformations of the form (\ref{transform}) with $c=0$ cannot change it. Next, analogously to (\ref{aux87}) we compute 
\begin{equation} 
\Theta^2_{21}=-\frac{S_{1}^{**}}{\sqrt{\rho_{11}^{**}}S^{**}}-\frac{\rho_{111}}{\sqrt{\rho_{11}^{3**}}},\label{aux78}
\end{equation}
and (\ref{finalxv}), (\ref{aux78}) yield
\begin{equation}
\begin{array}{l}
\displaystyle\tilde\Theta^2_{21}=\frac{1}{3S^{**}}\left[\frac{\rho_{12}^{**}}{\rho_{11}^{**}}\left(\frac{S_{1}^{**}}{\sqrt{\rho_{11}^{**}}S^{**}}\right)_{\hspace{-0.1cm}1}-\left(\frac{S_{1}^{**}}{\sqrt{\rho_{11}^{**}}S^{**}}\right)_{\hspace{-0.1cm}2}\right]-\\
\vspace{-0.1cm}\\
\displaystyle\hspace{4cm} \frac{1}{3S^{**}}\left[\frac{\rho_{12}^{**}}{\rho_{1 1}^{**}}\left(\frac{\rho_{111}}{\sqrt{\rho_{11}^{3**}}}\right)_{\hspace{-0.1cm}1}-\left(\frac{\rho_{111}}{\sqrt{\rho_{11}^{3**}}}\right)_{\hspace{-0.1cm}2}\right]-\frac{11S_{1}^{**}}{6\sqrt{\rho_{11}^{**}}\,\,{S^{**}}}-\frac{\rho_{111}}{6\sqrt{\rho_{11}^{3**}}}.
\end{array}\label{veryfinaltheta}
\end{equation}

It is now not hard to find a function $\rho$ for which the expression in the right-hand side of (\ref{veryfinaltheta}) does not identically vanish. Indeed, all solutions to the real homogeneous Monge-Amp\`ere equation can be explicitly described in parametric form (see, e.g. \cite{U}). For example, the solution corresponding to the choice $f(u)=u^3$, $g(u)=u$ in formula (2) of \cite{U} yields a function $\rho$ with the required properties. Indeed, in this case we have
$$
\rho=\frac{(1-12t_1t_2)^{\frac{3}{2}}+18t_1t_2-1}{108 t_2^2},\quad \rho_{11}=\frac{1}{\sqrt{1-12t_1t_2}},
$$
where $t_1$, $t_2$ are sufficiently small. Furthermore,
$$
S=\frac{1-\sqrt{1-12t_1t_2}}{t_2\sqrt{1-12t_1t_2}}\Biggl|_{t_{\alpha}=z_{\alpha}+{\bar z}_{\alpha}},
$$
hence $M$ is uniformly Levi degenerate of rank 1 in a neighborhood of the origin away from the set $\{\Re z_1=0\}$. In this case, computing the right-hand side of (\ref{veryfinaltheta}) we obtain
$$      
\tilde\Theta^2_{21}=-\frac{12 t_2}{(1-12t_1t_2)^{\frac{3}{4}}(1-\sqrt{1-12t_1t_2})}\Biggl|_{t_{\alpha}=z_{\alpha}+{\bar z}_{\alpha}},
$$
which shows that $\tilde\Theta^2_{21}$ does not identically vanish on $\gamma_0$.

\section{Applications of Theorem \ref{main}}\label{applications}
\setcounter{equation}{0}

In this section we discuss further properties of $\omega_M$ and give applications of our main result. 

First of all, by inspection of our construction in Section \ref{construction} one observes that for $M=\Gamma$ it leads to the bundle $G\ra H\hspace{-0.1cm}\setminus\hspace{-0.1cm}G$ and the Maurer-Cartan form $\omega_G^{\hbox{\tiny MC}}$. More precisely, upon identification of $\Gamma$ with $H\hspace{-0.1cm}\setminus\hspace{-0.1cm}G$, there exists an isomorphism $F$ of the bundles ${\mathcal P}_{\Gamma}\ra\Gamma$ and $G\ra H\hspace{-0.1cm}\setminus\hspace{-0.1cm}G$ that induced the identity map on the base and such that $F^*\omega_G^{\hbox{\tiny MC}}=\omega_{\Gamma}$. Hence the Maurer-Cartan equation for $\omega_G^{\hbox{\tiny MC}}$ yields $\Omega_{\Gamma}\equiv 0$.

Next, recall that $\pi_1(\Gamma)\simeq\ZZ_2$ (see p. 69 in \cite{FK1}) and consider the universal (double) cover of $\Gamma$. It can be realized as follows. First, note that the group $\hat G:=\Sp(4,\RR)$ is a double cover of $G$. Indeed, in order to see that $G$ is isomorphic to $\hat G/\{\pm\hbox{Id}\}$, one can realize the domains $\Omega_{\pm}$ introduced in (\ref{domainsomega}) in tube form (see \cite{FK1} and p. 289 in \cite{Sa}) and observe that the quotient $\hat G/\{\pm\hbox{Id}\}$ acts effectively on the tube realizations (see pp. 51--52 in \cite{FK1}). Now, let $\rho:\hat G\ra G$ be a 2-to-1 covering homomorphism and $\hat H:=\rho^{-1}(H)^{\circ}$. The subgroup $\hat H$ is isomorphic to $H$, and the induced map $\hat G/\hat H\ra G/H$ is a 2-to-1 covering. Thus, the quotient $\hat G/\hat H$ is a simply-connected covering space of $\Gamma$ endowed with an effective action of $\hat G$. We pull back the CR-structure from $\Gamma$ to this quotient and denote the resulting CR-manifold by $\hat\Gamma$. Clearly, one has $\Omega_{\hat\Gamma}\equiv0$. Further, analogously to (\ref{action}), define a right action of $\hat G$ on $\hat\Gamma$ by
\begin{equation}
\hat G\times\hat\Gamma\ra\hat\Gamma,\quad (g,p)\mapsto g^{-1}p,\label{hataction}
\end{equation}
and identify $\hat\Gamma$ with the right coset space $\hat H\hspace{-0.1cm}\setminus\hspace{-0.1cm}\hat G$ by means of this action. For the manifold $\hat\Gamma$ our construction in Section \ref{construction} leads to the bundle $\hat G\ra \hat H\hspace{-0.1cm}\setminus\hspace{-0.1cm}\hat G$ and the right-invariant Maurer-Cartan form $\omega_{\hat G}^{\hbox{\tiny MC}}$.

We call a manifold $M\in{\mathfrak C}_{2,1}$ flat if $\Omega_M\equiv 0$. For instance, both $\Gamma$ and $\hat\Gamma$ are flat. Further, we say that $M$ is locally CR-equivalent to $\Gamma$ if for every point $p\in M$ there exists a neighborhood of $p$ that is CR-equivalent to an open subset of $\Gamma$. In our first corollary to Theorem \ref{main} we show that these classes of manifolds in fact coincide.

\begin{corollary}\label{cor1}\sl A manifold $M\in{\mathfrak C}_{2,1}$ is flat if and only if $M$ is locally CR-equivalent to $\Gamma$.
\end{corollary}

\noindent{\bf Proof:} If $\Omega_M\equiv 0$, then for any ${\bf\Theta}\in{\mathcal P}_M$ there exists a diffeomorphism $F:U\ra V$, where $U$ is a neighborhood of ${\bf\Theta}$ in ${\mathcal P}_M$ and $V$ is a neighborhood of the identity in $G$, such that $F^*\omega_G^{\hbox{\tiny MC}}=\omega_M$ (see, e.g. Theorem 1.2.4 in \cite{CSl}). By Theorem \ref{main}, the map $F$ is a lift of a CR-isomorphism between a neighborhood of $\pi({\bf\Theta})\in M$ and an open subset of $H\hspace{-0.1cm}\setminus\hspace{-0.1cm}G\simeq\Gamma$, and therefore $M$ is locally CR-equivalent to $\Gamma$.

Conversely, suppose that $M$ is locally CR-equivalent to $\Gamma$. Then Theorem \ref{main} and the flatness of $\Gamma$ imply that $\Omega_M\equiv 0$. \qed 
\vspace{0.3cm}

Our next corollary concerns the extendability of local CR-automorphisms of $\Gamma$ and $\hat\Gamma$. As we mentioned in Section \ref{model}, this result is known and follows from Theorems 4.5, 4.7 in \cite{KZ}. Here we give a new short proof based on Theorem \ref{main}.

\begin{corollary}\label{cor2}\sl Any local CR-auto\-mor\-phism of either $\Gamma$ or $\hat\Gamma$ extends to a global CR-auto\-mor\-phism induced by the action of either $G$ or $\hat G$, respectively. In particular $\Aut(\Gamma)\simeq G$ and $\Aut(\hat\Gamma)\simeq\hat G$.
\end{corollary}

\noindent {\bf Proof:} The proof is analogous to that of Theorem 6 in \cite{T1} for Levi nondegenerate hyperquadrics. Let ${\bf\Gamma}$ be either $\Gamma$ or $\hat\Gamma$ and ${\bf G}$ be either $G$ or $\hat G$, respectively. Consider a CR-isomorphism $f:U\ra V$, where $U,V\subset{\bf \Gamma}$ are domains. By Theorem \ref{main}, the map $f$ can be lifted to a bundle isomorphism $F:\pi^{-1}(U)\ra\pi^{-1}(V)$ preserving the absolute parallelism $\omega_{{\bf\Gamma}}$. Identifying $\pi^{-1}(U)$, $\pi^{-1}(V)$ with domains in ${\bf G}$ and $\omega_{{\bf \Gamma}}$ with $\omega_{\bf G}^{\hbox{\tiny MC}}$, we then see that $F$ is given by the right multiplication by an element $g\in{\bf G}$. Now (\ref{action}), (\ref{hataction}) imply that $f$ is induced by the (left) action of $g^{-1}$ on ${\bf\Gamma}$.\qed
\vspace{0.3cm}

The next corollary concerns germs of CR-isomorphisms between germs of arbitrary manifolds in ${\mathfrak C}_{2,1}$. For $M\in{\mathfrak C}_{2,1}$ and $p\in M$ the germ of $M$ at $p$ is denoted by $(M,p)$.

\begin{corollary}\label{cor4}\sl Let $f_1,f_2:(M,p)\ra(\tilde M,\tilde p)$ be CR-isomorphism germs. Assume that the 2-jets $j^2(f_1)$ and $j^2(f_2)$ of $f_1$ and $f_2$ coincide. Then $f_1=f_2$.
\end{corollary}

\noindent{\bf Proof:} For $k=1,2$ let $\hat f_k: U_k\ra \tilde U_k$ be a CR-isomorphism representing the germ $f_k$, where $U_k$ is a neighborhood of $p$ in $M$ and $\tilde U_k$ a neighborhood of $\tilde p$ in $\tilde M$. Further, let $F_k:{\mathcal P}_{U_{{}_k}}\ra{\mathcal P}_{\tilde U_{{}_k}}$ be the lift of $\hat f_k$ such that $F_k^*\omega_{\tilde U_{{}_k}}=\omega_{U_{{}_k}}$. Since  $j^2(\hat f_1,p)=j^2(\hat f_2,p)$, it follows from our construction in Section \ref{construction} that $F_1\equiv F_2$ on $\pi^{-1}(p)$. This implies that in fact $F_1\equiv F_2$ in a neighborhood of $\pi^{-1}(p)$ in ${\mathcal P}_M$ (see Lemma 5.14 in \cite{E}). Hence $f_1=f_2$.\qed
\vspace{-0.3cm}\\

\noindent As the example of the manifold $\Gamma$ shows, Corollary \ref{cor4} is sharp in the sense that the 2-jet determination property obtained there cannot be replaced by the 1-jet determination property.

Further, for $M\in{\mathfrak C}_{2,1}$ and $p\in M$, let $\Stab(M,p)$ be the stability group of $M$ at $p$, i.e.~the group of all CR-isomorphism germs from $(M,p)$ to itself. We endow $\Stab(M,p)$ with the topology of 2-jet evaluation at the point $p$, which turns $\Stab(M,p)$ into a topological group. Also, let $\aut(M,p)$ be the Lie algebra of germs of infinitesimal CR-automorphisms of $M$ at $p$. In the next corollary we show that any stability group can be viewed as a subgroup of $H$ and give an estimate on $\dim\aut(M,p)$.  

\begin{corollary}\label{cor5}\sl For any $M\in{\mathfrak C}_{2,1}$ and $p\in M$ there exists a continuous injective homomorphism ${\bf F}_p:\Stab(M,p)\ra H$, and $\dim\aut(M,p)\le 10$. 
\end{corollary}

\noindent{\bf Proof:} For $f\in\Stab(M,p)$ let $\hat f: U\ra V$ be a CR-isomorphism representing the germ $f$, where $U$ and $V$ are neighborhoods of $p$ in $M$. Further, let $F:{\mathcal P}_U\ra{\mathcal P}_V$ be the lift of $\hat f$ such that $F^*\omega_V=\omega_U$. Clearly, $F$ preserves the fiber $\pi^{-1}(p)$. We now restrict $F$ to $\pi^{-1}(p)$ and denote the restriction by the same symbol. The group $H$ acts on the fiber freely transitively, and we identify $\pi^{-1}(p)$ with $H$ by means of this action. Upon this identification, it follows from property (\ref{ver}) that $F$ preserves all right-invariant vector fields on $H$ and therefore is given by the right multiplication by an element $h\in H$. We then define ${\bf F}_p$ as the map $f\mapsto h^{-1}$. It is clear that ${\bf F}_p$ is a continuous homomorphism, and the injectivity of ${\bf F}_p$ follows by the argument given at the end of the proof of Corollary \ref{cor4}. The estimate for $\dim\aut(M,p)$ is obtained as in Theorem 2.6 in \cite{BS2}. \qed
\vspace{-0.3cm}\\

\noindent Corollary \ref{cor5} together with Theorem VII of \cite{P} implies that $\Stab(M,p)$ admits the structure of a Lie group with Lie algebra $\{v\in{\mathfrak h}:\exp(tv)\in{\bf F}_p(\Stab(M,p))\,\hbox{for all $t\in\RR$}\}$ and possibly uncountably many connected components. It is not clear, however, whether the Lie group topology always coincides with the topology of 2-jet evaluation at $p$ introduced on $\Stab(M,p)$ earlier.

Our next result focusses on the group $\Aut(M)$ of global CR-automorphisms of a manifold $M\in{\mathfrak C}_{2,1}$, namely on the existence of a Lie group structure on $\Aut(M)$. In what follows $\Aut(M)_p$ is the isotropy subgroup of a point $p\in M$ under the $\Aut(M)$-action and ${\bf F}_p':={\bf F}_p\circ\iota$, where ${\bf F}_p$ is the homomorphism constructed in the proof of Corollary \ref{cor5} and\linebreak $\iota:\Aut(M)_p\ra\Stab(M,p)$ is the homomorphism that assigns to every element of $\Aut(M)_p$ its germ at $p$. 

\begin{corollary}\label{cor3}\sl For $M\in{\mathfrak C}_{2,1}$ the group $\Aut(M)$ admits the structure of a Lie transformation group of $M$ of dimension at most 10 such that the following holds:

\noindent{\rm (i)} the Lie algebra of $\Aut(M)$ is isomorphic to the Lie algebra $\aut_c(M)$ of complete infinitesimal CR-automorphisms of $M$;
\vspace{0.1cm}\\
\noindent{\rm (ii)} for any $p\in M$ the map ${\bf F}_p'$ is a Lie group homomorphism with respect to the induced Lie group structure on $\Aut(M)_p$;
\vspace{0.1cm}\\
\noindent{\rm (iii)} the Lie group topologies on $\Aut(\Gamma)$ and $\Aut(\hat\Gamma)$ coincide with the compact-open topologies and also with the topologies arising from the Lie groups $G$, $\hat G$, respectively;
\vspace{0.1cm}\\
\noindent{\rm (iv)} if $\dim\Aut(M)=10$, then $M$ is CR-equivalent to either $\Gamma$ or $\hat\Gamma$.
\end{corollary}

\noindent{\bf Proof:} For $M\in{\mathfrak C}_{2,1}$, let $\Diff({\mathcal P}_M)$ be the group of diffeomorphisms of the fiber bundle ${\mathcal P}_M$ endowed with the compact-open topology. By Theorem \ref{main}, every $f\in\Aut(M)$ lifts to a uniquely defined element $F\in\Diff({\mathcal P}_M)$ that preserves the absolute parallelism $\omega_M$, and it is clear that the map $f\mapsto F$ is a homomorphism. Hence $\Aut(M)$ is isomorphic as an abstract group to the subgroup $\Diff({\mathcal P}_M)_{\omega_{{}_M}}$ of $\Diff({\mathcal P}_M)$ that consists of all diffeomorphisms of $M$ preserving $\omega_M$. By a result due to S. Kobayashi (see Theorem 3.2 in \cite{Kob}), the group $\Diff({\mathcal P}_M)_{\omega_{{}_M}}$ acts freely on ${\mathcal P}_M$ and admits the structure of a Lie group of dimension at most 10 with respect to the compact-open topology whose Lie algebra is isomorphic to $\aut_c(M)$. Transferring this Lie group structure to $\Aut(M)$, one observes that $\Aut(M)$ acts on $M$ as a topological group. Hence with respect to this structure $\Aut(M)$ is a Lie transformation group of $M$ (see pp. 212--213 in \cite{MZ}). Therefore, $\Aut(M)_p$ is a closed subgroup of $\Aut(M)$, and the Lie group structure on $\Aut(M)$ induces a Lie group structure on $\Aut(M)_p$. It is clear from the above description of the topology on $\Aut(M)$ that for any $p\in M$ the map ${\bf F}_p'$ is continuous and thus is a Lie group homomorphism. It also follows that the Lie group topologies constructed on $\Aut(\Gamma)$ and $\Aut(\hat\Gamma)$ as above coincide with the corresponding compact-open topologies and also with the topologies arising from the isomorphisms $\Aut(\Gamma)\simeq G$, $\Aut(\hat\Gamma)\simeq\hat G$.   

Next, let $\dim\Aut(M)=10$. This condition means that $\Diff({\mathcal P}_M)_{\omega_{{}_M}}$ acts freely transitively on ${\mathcal P}_M$. Hence $\Aut(M)$ is a connected group that acts on $M$ transitively, and we can view $M$ as a real-analytic manifold for which the 10-dimensional algebra $\aut_c(M)$ consists of real-analytic vector fields. Theorems I and II in \cite{FK2} now imply that $M$ is locally CR-equivalent to $\Gamma$. Let $\hat M$ be the universal cover of $M$. Clearly, one has $\dim\Aut(\hat M)=10$. As in the case of Levi nondegenerate CR-hypersurfaces locally CR-equivalent to the sphere, one can define a development map ${\mathcal F}:\hat M\ra\Gamma$ and a development representation $\Pi:\Aut(\hat M)\ra\Aut(\Gamma)$ (see pp. 224--225 in \cite{BS1}), where for any $g\in\Aut(\hat M)$ we have
$$
{\mathcal F}\circ g=\Pi(g)\circ{\mathcal F}.
$$
Since $\hat M$ is homogeneous, arguing as in the proof of Proposition 1.4 in \cite{BS1}, one obtains that ${\mathcal F}$ is a covering map from $\hat M$ onto ${\mathcal F}(\hat M)$ (cf. Proposition 6.4 in \cite{FK1}). Furthermore, the homomorphism $\Pi$ is continuous and has discrete kernel. Therefore, $\Pi$ is surjective, which implies ${\mathcal F}(\hat M)=\Gamma$. It then follows that $\hat M$ is CR-equivalent to $\hat\Gamma$. Further, we have $\Aut(\hat\Gamma)\simeq\hat G$, hence the group of deck transformations of the covering $\hat\Gamma\ra M$ is isomorphic to a central subgroup of $\hat G$. Since the center of $\hat G$ is $\pm\hbox{Id}$, it follows that $M$ is CR-equivalent to either $\hat\Gamma$ or $\Gamma$.\qed
\vspace{-0.3cm}\\

For an arbitrary smooth CR-manifold that is finitely nondegenerate and of finite type in the sense of Kohn and Bloom-Graham the existence of a Lie group structure on $\Aut(M)$ with respect to the $C^{\infty}$ compact-open topology was established in Theorem 6.2 of \cite{BRWZ}. Observe that in Corollary \ref{cor3} the topology on $\Aut(M)$ is a priori weaker than the $C^{\infty}$ compact-open topology. We also remark that the Lie group structure on $\Aut(M)$ constructed in Corollary \ref{cor3} can be obtained by arguing as in Corollary 2.8 in \cite{BS2} and using Corollary \ref{cor5}.

Finally, we will settle a conjecture due to V. Beloshapka for 5-dimensional manifolds. For every $n\ge 2$ let ${\mathfrak C}_n$ be the class of all connected real-analytic CR-hypersurfaces $M$ of dimension $2n-1$ with $\dim\hol(M,p)<\infty$ for all $p\in M$, where $\hol(M,p)\subset\aut(M,p)$ is the Lie subalgebra of germs of real-analytic infinitesimal CR-automorphisms of $M$ at $p$. Note that\linebreak $\hol(M,p)=\aut(M,p)$ if, for instance, $M$ is essentially finite at $p$ (see Theorem 2.2 and Remark 2.5 in \cite{Sta}).   

\begin{conjecture}\label{belo} {\rm (V. Beloshapka)} \sl For any manifold $M\in{\mathfrak C}_n$ and any $p\in M$ one has\linebreak $\dim\hol(M,p)\le n^2+2n$, and the maximal possible value $n^2+2n$ of $\dim\hol(M,p)$ is attained at some point of $M$ if and only if $M$ is generically locally $CR$-equivalent to a Levi nondegenerate hyperquadric in $\CC^n$.
\end{conjecture} 

As examples in \cite{Kow}, \cite{Z} show, a manifold $M\in{\mathfrak C}_n$ with $\dim\hol(M,p_0)=n^2+2n$ at some point $p_0\in M$ need not be locally CR-equivalent to a hyperquadric everywhere. On the other hand, we make the following remark.

\begin{remark}\label{conjremark}\rm If $M$ is Levi nondegenerate at $p_0$ and $\dim\hol(M,p_0)\ge n^2+2n$, then\linebreak $\dim\hol(M,p_0)=n^2+2n$ and $M$ is CR-equivalent to the corresponding hyperquadric $Q$ near $p_0$. Indeed, results of \cite{CM} imply that if the Levi form of $M$ is nondegenerate at a point $q$, then $\Stab(M,q)$ can be embedded in the Lie group $\Stab(Q,0)$ as a closed subgroup with Lie algebra isomorphic to the subalgebra $\hol_0(M,q)\subset\hol(M,q)$ of all real-analytic infinitesimal CR-automorphism germs vanishing at $q$. A detailed analysis of the embedding shows that if $\dim\Stab(M,q)\ge n^2-2n+2$ (which is weaker than the condition $\dim\hol(M,q)\ge n^2+2n$) then $M$ is CR-equivalent to $Q$ in a neighborhood of $q$ (see Section 9 in \cite{V}).
\end{remark}   

For $n=2$ the conjecture follows from results of \cite{C1} since for any $M\in{\mathfrak C}_2$ the set of points where $M$ is Levi nondegenerate is nonempty and therefore dense in $M$. If $n\ge 3$ a manifold $M\in{\mathfrak C}_n$ may not have any points of Levi nondegeneracy, thus the above argument does not always work. Below we confirm the conjecture for $n=3$. In the special case when $M\in{\mathfrak C}_3$ is locally homogeneous, the conjecture also follows from results of \cite{FK2}. To the best of our knowledge, it remains open for $n>3$.  

\begin{corollary}\label{cor8}\sl Conjecture {\rm \ref{belo}} holds for $n=3$.
\end{corollary}

\noindent {\bf Proof:} Let $M\in{\mathfrak C}_3$ and suppose that there exists $p_0\in M$ with $\dim\hol(M,p_0)\ge 15$. Denote by $S$ the set of points at which $M$ is Levi nondegenerate. If $S=\emptyset$, then either the manifold $M$ is Levi-flat or there exists a proper real-analytic subset $V\subset M$ such that every connected component of $M\setminus V$ lies in ${\mathfrak C}_{2,1}$. In the former case $\dim\hol(M,p)=\infty$ for all $p\in M$, which is impossible. In the latter case by choosing a point $p\in M\setminus V$ sufficiently close to $p_0$ we obtain $\dim\hol(M,p)\ge 15$, which contradicts the last statement of Corollary \ref{cor5}. 

Thus, we have shown that $S\ne\emptyset$, hence $S$ is dense in $M$. Therefore, there exists a point $p\in S$ arbitrarily close to $p_0$ such that $\dim\hol(M,p)\ge \dim\hol(M,p_0)$. Application of Remark \ref{conjremark} with $n=3$ now yields that $\dim\hol(M,p)= \dim\hol(M,p_0)=15$ and that $M$ is CR-equivalent to the corresponding hyperquadric in $\CC^3$ near the point $p$.\qed
\vspace{0.3cm}

{\bf Acknowledgements.} We are grateful to A. \v Cap for many useful conversations. The second author would like to thank the Australian National University for its hospitality during his visit to Canberra in 2012 when this work was initiated.

{\obeylines
\noindent A. Isaev: Department of Mathematics, The Australian National University, Canberra, 

\noindent ACT 0200, Australia, e-mail: {\tt alexander.isaev@anu.edu.au}
\vspace{0.3cm}

\noindent D. Zaitsev: School of Mathematics, Trinity College Dublin, Dublin 2, Ireland,

\noindent e-mail: {\tt zaitsev@maths.tcd.ie}

}

\end{document}